\documentclass{amsart}
\begin{document}
\allowdisplaybreaks
\newtheorem{theorem}{{\bf Theorem}}
\newtheorem{problem}{{\bf Problem}}
\newtheorem{lemma}{{\bf Lemma}}
\newtheorem{cor}{{\bf Corollary}}
\newtheorem{definition}{{\bf Definition}}
\makeatletter
\def\labelenumi{(\@roman\c@enumi)}
\def\theenumi{(\@roman\c@enumi)}
\def\labelenumii{(\@alph\c@enumii)}
\def\theenumii{\@alph\c@enumii}
\def\labelenumiii{(\@arabic\c@enumiii)}
\def\theenumiii{(\@arabic\c@enumiii)}
\def\p@enumiii{\theenumi\theenumii}
\def\wideA#1{\@mathmeasure\z@\textstyle{#1}\ifdim\wd\z@>\tw@ em\mathaccent "055C{#1}%
\else\mathaccent "0364{#1}\fi}
\def\wideB#1{\@mathmeasure\z@\textstyle{#1}\ifdim\wd\z@>\tw@ em\mathaccent "055E{#1}%
\else\mathaccent "0366{#1}\fi}
\makeatother
\def \blank{\phantom{x}}
\def \con#1{\setbox13\hbox{$#1$}\ifdim\wd13<1em\breve{#1}\else{\(#1\)}\breve{\ }\fi}
\def \concept#1{{\bf#1}}
\def \halfthinspace{\relax\ifmmode\mskip.5\thinmuskip\relax\else\kern.8888em\fi}
\let \hts=\halfthinspace
\def \rp{{\hts;\hts}}
\let \bp=\cdot
\def \id{{1\kern-.08em\raise1.3ex\hbox{\rm,}\kern.08em}}
\def \di{{0\kern-.04em\raise1.3ex\hbox{\rm,}\kern.04em}}
\let \SS=\S
\def \makecs#1#2{\makecsX {#1}#2,.}
\def \makecsX#1#2#3.{\onecs{#1}{#2} 
\ifx#3,\let\next\eatit\else\let\next\makecsX\fi\next{#1}#3.}
\def \onecs#1#2{\expandafter\gdef\csname #2\endcsname%
{{\csname #1\endcsname {#2}}}}
\def \eatit#1#2.{\relax}
\makecs{}{abcdefghijklmnopqrstuvwxyzABCDEFGHIJKLMNOPQRSTUVWXYZ}
\def \gc#1{\mathbf{#1}}   
\def \Sb#1{Sb\hts\(#1\)}
\def \comp#1{\overline{#1}}
\def \Fd#1{\setbox13\hbox{$#1$}\ifdim\wd13<1pt{Fd\hts#1}\else{Fd\(#1\)}\fi}
\def \Do#1{\setbox13\hbox{$#1$}\ifdim\wd13<1pt{Do\hts#1}\else{Do\(#1\)}\fi}
\def \Ra#1{\setbox13\hbox{$#1$}\ifdim\wd13<1pt{Ra\hts#1}\else{\gc{Ra}\(#1\)}\fi}
\def \Ca#1{\gc{Ca}\(#1\)}
\def \Na#1{\gc{Nr}_{#1}}
\def \({\left(}
\def \){\right)}
\def \<{\left<}
\def \>{\right>}
\def \minus{\mathop\sim}
\def \min#1{\overline{#1}}
\def \RRA{\mathsf{RRA}}
\def \wRRA{\mathsf\w\RRA}
\def \At#1{\gc{At}\(#1\)}
\def \at#1{At\(#1\)}
\def \RA{\mathsf{RA}}
\def\CA#1{\mathsf{CA}_{#1}}
\def \nr{{\mathsf{Nr}}}
\def \Nr{{\gc{Nr}}}
\def\Nrr#1{\mathsf{Nr}_{#1}}
\def \Sub{\mathbf\S}
\def \sub#1#2{{\sf s}^{#1}_{#2}}
\def \cyl#1{{\mathsf c}_{#1}}
\def \diag#1#2{{\mathsf d}_{#1#2}}
\def \etc{{\it etc}}
\def \NA{\mathsf{NA}}
\def \SA{\mathsf{SA}}
\def \Cm#1{\setbox13\hbox{$#1$}\ifdim\wd13=0pt{\gc{Cm}}\else{\gc{Cm}\({#1}\)}\fi}
\def \CmB{\gc{Cm\,At(\B)}}
\def\Sg#1#2{\gc{Sg}^{(#1)}\hts({#2})}
\def\nul{\text{---}}
\def\lll{\hline\rule[10pt]{0pt}{0pt}}
\def\alg#1#2{\mathsf{[#1]}\mathsf{#2}}
\def\alg#1#2{\text{$\mathsf{#1}_{#2}$}}
\def\atsymbol{\char'100}
\def\?{\cdot}
\def\mo{\mathbf\mu}
\def\dd#1{\mathsf{d}}
\def\aa#1{\mathsf{a}_{#1}}
\def\cc#1{\mathbf{e}_{#1}}
\def\er{\gc\E_\r^{123}}
\def\ep{\gc\E^{23}_{\r+3}}
\def\eq{\gc\E^{23}_\q}
\def\ex{\gc\E^{23}_7}
\def\cov#1{\setbox13\hbox{$#1$}\ifdim\wd13<.1em{\mathsf{c}{#1}}\else{\mathsf{c}(#1)}\fi}
\def\Te#1#2#3{\T(#1,#2,#3)}
\def\bn{\B_\n}
\def\bo{\A\t_0}
\def\join{\sum}
\def\ie{{\it i.e.\/}}
\def\jn#1#2{\J({#1},{#2})}
\def\an{\jn\a\n}
\def\ao{\jn\a0}
\def\CC{\gc\C}
\def\RR{\gc\B}
\def\pr#1#2{#1^{(#2)}}
\def\pra#1#2{\mathsf{e}_{#1}^{(#2)}}
\def\bu{\bullet}
\def \algN{\gc\B}
\def\sq{S^w_q}
\def\hc#1{\mathbf{u}_{#1}}
\def\hc#1{u^h_#1}
\def\hc#1{t_{#1}}
\def\hyph{-}
\def\lb{\left\{\,}
\def\rb{\,\right\}}
\def\gSb#1{\gc{Sb}\hts\(#1\)}
\def\univ{\mathsf\V}
\def\al{\alpha}
\def\ka{\kappa}
\def\la{\lambda}
\def\ls#1#2{L\(#1,#2\)}
\def\ee#1{1_{#1}}
\def\epp#1{e_{#1}}
\def\dr#1#2#3{d^{#1}_{#2#3}}
\def\Dr#1#2#3{D^{#1}_{#2#3}}
\def\ddd#1#2#3#4{d^{#1(#2)}_{#3#4}}
\def\dash{{\text{-}}}
\def\CPQ{\T^*(\P,\Q)}
\def\cy#1{c_{#1}}
\def\su#1#2{s^{#1}_{#2}}
\def\gen{G_0}
\def\ratios{{\mathbb Q}}
\def\CEA{\C_\gc{E}(\gc{A})}
\def\CerA{\C_{\er}\(\gc\A\)}
\def\Cerep{\C_{\er}\(\ep\)}
\def\false{\sf false}
\def\true{\sf true}
\def\MM{\mathcal{M}}
\def\mod{\mathcal{M}}
\def\nog{\text{not }}
\def\tsum{\textstyle\sum}
\def\Cex{\C_\gc{E}(\ex)}
\def\uu#1{u_{#1}}
\def\vv#1{v_{#1}}
\def\As{\gc\A^\omega}
\def\lang#1{\L^{\gc\A}_{#1}}
\def\thA{\T\h(\gc\A)}

\title[Subcompletions]{Subcompletions of representable relation algebras}
\author{Roger D. Maddux}
\address{Department of Mathematics
\\396 Carver Hall 
\\Iowa State University 
\\Ames, Iowa 50011-2066, USA}
\email{ maddux@iastate.edu }
\thanks{}
\date{August 3, 2011}
\begin{abstract}
  Many finite symmetric integral non-representable relation algebras, including almost
  all Monk algebras, can be embedded in the completion of an atomic symmetric integral
  representable relation algebra whose finitely-generated subalgebras are finite.
\end{abstract}

\maketitle

\section{\bf Introduction}
Monk~\cite{MR0277369} proved that if $\gc\B$ is a Boolean algebra with operators, then
$\gc\B$ has a unique completion $\CC$, where $\CC$ is a {\bf completion} of $\gc\B$ if
$\CC$ is a complete, $\gc\B$ is a subalgebra of $\CC$, and $\gc\B$ is \concept{dense} in
$\CC$, which means that below every non-zero element of $\CC$ there is a non-zero
element of $\gc\B$. Monk proved that if $\gc\B$ is a relation algebra algebra, then its
completion is also relation algebra. The problem remained, if $\gc\B$ is representable,
must its completion also be representable? In other words, is the variety\footnote{It is
easy to show that $\RRA$ is closed under subalgebras and direct products. Closure under
homomorphic images was first proved by Tarski~\cite{MR0066303} using model theory; for a
more direct proof see~\cite[Th.\,121]{MR2269199}.} $\RRA$ of representable relation
algebras closed under completions?  Hodkinson~\cite{MR1490103} provided the answer that,
no, $\RRA$ is not closed under completions because there is an \emph{atomic}
$\gc\B\in\RRA$ such that the completion of $\gc\B$ is not representable.\footnote{For
subsequent developments, alternate and simpler constructions, extensions, and related
work, see, for example, Andr\'eka-N\'emeti-Sayed~Ahmed~\cite{ MR2387933},
Hirsch-Hodkinson~\cite{MR1330986, MR1472125, MR1935083, MR1887031, MR2548463},
Hodkinson-Venema~\cite{ MR2156722}, Khaled-Sayed~Ahmed~\cite{MR2605423, MR2507303},
Sayed~Ahmed~\cite{MR2357188, MR2443834, MR2441102, MR2417802, MR2519240}, and
Sayed~Ahmed-Samir~\cite{ MR2453362}.}

Consider any atomic representable relation algebra $\gc\B\in\RRA$ whose completion $\CC$
is not representable, as might arise from Hodkinson's proof. Since $\RRA$ is a variety
and $\CC$ is not in $\RRA$, there must be an equation $\epsilon$ that holds in $\RRA$
but fails in $\CC$.  Let $\gc\A$ be the subalgebra of $\CC$ that is generated by the
finitely many values assigned to the variables occurring in $\epsilon$.  Then $\gc\A$ is
also not representable because it fails to satisfy $\epsilon$.  Thus $\gc\A$ is an
example of a finitely-generated relation algebra which is a subalgebra of the
non-representable completion of an atomic $\RRA$. The question addressed by this paper
is, which relation algebras can occur as $\gc\A$? We rephrase this question as a
problem.
\begin{problem}\label{prob}
  Let $\K$ be the class of finitely-generated subalgebras of non-representable
  completions of atomic representable relation algebras. Which relation algebras are in
  $\K$? Does $\K$ contain any relation algebras that are not weakly representable?
\end{problem}
Some finite symmetric integral relation algebras have no proper extensions at all and
are therefore neither representable nor in $\K$; see Frias-Maddux~\cite{MR1608984} for
examples.

As a partial positive answer, it will be shown in this paper that every finite Monk
algebra with six or more colors is in $\K$. Monk algebras are relation algebraic
versions of cylindric algebras used by Monk~\cite{MR0256861} to prove that classes of
finite-dimensional representable cylindric algebras are not finitely axiomatizable.

\section{\bf Monk algebras}
\begin{definition}\label{def1}
  For $4\leq\q\in\omega$, the \concept{no 1-cycles algebra} $\eq$ is the finite
  symmetric integral relation algebra with $\q$ atoms $\cc0=\id$, $\cc1$, $\cdots$,
  $\cc{\q-1}$ such that if $\a,\b$ are distinct diversity atoms then $\a\rp\b=\di$ and
  $\a\rp\a=\comp\a$.
\end{definition}
These relation algebras were constructed in~\cite{Maddux1978} and were called
$\mathfrak{E}_\q\(\{2,3\}\)$ in~\cite[Def.\,2.4, Prob.\,2.7]{MR1011183}.  Two particular
examples of these algebras, namely $\gc\E^{23}_4$ and $\gc\E^{23}_5$, are \alg{62}{65}
and~\alg{3009}{3013}, respectively, in~\cite{MR2269199}.  It is likely that $\eq$ is
representable for all $\q\geq4$. This is known to be true for $\q=4$; the earliest
reference to the representability on a 13-element set of $\gc\E^{23}_4$ is in a footnote
by Lyndon~\cite{MR0037278}.  In fact, $\gc\E^{23}_4$ is isomorphic to a subalgebra of
the complex algebra of the 13-element cycle group $\mathbb{Z}_{13}$, and $\gc\E^{23}_4$
has exactly two other square representations, both on a 16-element set. Furthermore,
$\eq$ was shown to be representable by Comer~\cite{MR734546} for $\q=5,6$ and more
recently (by Comer's method) for $\q=7,8$.

$\eq$ can also be described by cycles and atom structures, which are defined for
algebras in $\NA$, the class of~\concept{non-associative relation algebras}.  An axiom
set for $\NA$ is obtained by deleting the associative law from the axioms for relation
algebras; see~\cite[Def.\,1.2]{MR662049} and~\cite[Th.\,314]{MR2269199}. The~\concept{atom
structure}~\cite[Def.\,3.2]{MR662049} of an algebra $\gc\A\in\NA$ is
$\<\at{\gc\A},\C,\con\blank,\I\>$ where $\at{\gc\A}$ is the set of atoms of $\gc\A$,
$\C$ is the set of triples of atoms $\<\x,\y,\z\>$ such that $\x\rp\y\geq\z$,
$\con\blank$ is the restriction of the converse operation of $\gc\A$ to the atoms of
$\gc\A$, and $\I=\{\x:\id\geq\x\in\at{\gc\A}\}$.  In every $\NA$, $\C$ is the union of
sets of the form
\begin{equation}
[\x,\y,\z]=\{	\langle    \x,    \y,    \z\rangle,
		\langle\con\x,    \z,    \y\rangle,
		\langle    \y,\con\z,\con\x\rangle,
		\langle\con\y,\con\x,\con\z\rangle,
		\langle\con\z,    \x,\con\y\rangle,
		\langle    \z,\con\y,    \x\rangle\},
\end{equation}
where $\x,\y,\z\in\at{\gc\A}$. Such sets are called \concept{cycles}.  If $\id$ is an
atom of $\gc\A$, then the cycle $[\x,\y,\z]$ is said to be an \concept{identity cycle}
if $\id\in[\x,\y,\z]$, and a \concept{diversity cycle} otherwise. If $\gc\A$ is
\concept{symmetric}, \ie, $\con\x=\x$ for all $\x$, then a diversity cycle $[\x,\y,\z]$
is said to be a~\concept{1-cycle}, \concept{2-cycle}, or~\concept{3-cycle} if the
cardinality $|\{\x,\y,\z\}|$ is 1, 2, or 3, respectively.  For example, the cycles of
$\eq$ are all the 2-cycles and 3-cycles, but none of the 1-cycles.
\begin{definition}[{Andr\'eka-Maddux-N\'emeti~\cite{MR1052567}}]\label{def2}
  Let $\gc\A$ and $\gc\B$ be atomic relation algebras. We say that $\gc\A$ is obtained
  from $\gc\B$ by \concept{splitting} if $\gc\B\subseteq\gc\A$, every atom $\x$ of
  $\gc\A$ is contained in an atom $\cov\x$ of $\gc\B$, called the \concept{cover} of
  $\x$, and for all $\x,\y\in At\,\gc\A$, if $\x,\y\leq\di$ then
\begin{equation}\label{def2.1}
\x\rp\y=
\begin{cases} 
	\cov\x\rp\cov\y\bp\di&\text{if }\x\neq\con\y
\\	\cov\x\rp\cov\y &\text{if }\x=\con\y.
\end{cases}
\end{equation}
\end{definition}

\begin{definition}[{Andr\'eka-Maddux-N\'emeti~\cite[Ex.\,6]{MR1052567}}]\label{def3}
  A~\concept{Monk algebra} is an atomic symmetric integral relation algebra obtained by
  splitting from some $\eq$, $4\leq\q\in\omega$.
\end{definition}
Assume $\gc\A$ is a Monk algebra obtained from $\eq$ by splitting.  Then $\gc\A$ extends
$\eq$ and the $\q-1$ diversity atoms of the subalgebra $\eq\subseteq\gc\A$ are called
the~\concept{colors} of $\gc\A$. Consider a subalgebra
$\gc\E\subseteq\eq\subseteq\gc\A$.  Then the Monk algebra $\gc\A$ extends its subalgebra
$\gc\E$ in a way that, for want of a better name, we simply call ``special''.
\begin{definition}\label{spec}
  If $\gc\A$ and $\gc\E$ are finite symmetric integral relation algebras, then $\gc\A$
  is said to be a~\concept{special extension} of $\gc\E$ if $\gc\E\subseteq\gc\A$ and
  for all diversity atoms $\di\geq\a,\b,\c\in\at{\gc\E}$,
\begin{enumerate}
\item	if not $(\a=\b=\c)$ and $\a\rp\b\geq\c$ then $\x\rp\y\geq\c$ whenever
	$\a\geq\x\in\At{\gc\A}$ and $\b\geq\y\in\At{\gc\A}$,\label{spec.i}
\item	if $\a\rp\a\geq\a$ then $\x\rp\y\bp\a\neq0$ whenever
	$\a\geq\x,\y\in\At{\gc\A}$.\label{spec.ii}
\end{enumerate}
\end{definition}
Every finite symmetric integral relation algebra is a special extension of itself.
Every finite symmetric integral relation algebra with no functional atoms\footnote{A
functional atom $\x$ in a symmetric integral relation algebra satisfies $\x\rp\x=\id$.}
is also a special extension of its minimum subalgebra, the one whose atoms are $\id$ and
$\di$.
\begin{lemma}\label{lem1}
  Every Monk algebra obtained from $\eq$ by splitting is a special extension of every
  subalgebra of $\eq$.
\end{lemma}
\proof
Assume $\eq\subseteq\gc\A$, $4\leq\q$, $\gc\A$ is a Monk algebra obtained from $\eq$ by
splitting, and $\cov\x$ is the atom of $\eq$ containing the atom $\x$ of $\gc\A$.
Consider a subalgebra $\gc\E\subseteq\eq\subseteq\gc\A$ and diversity atoms
$\di\geq\a,\b,\c\in\at{\gc\E}$.

To show part~(i) of Def.\,\ref{spec}, we assume $\nog(\a=\b=\c)$, $\a\rp\b\geq\c$,
$\a\geq\x\in\at{\gc\A}$, and $\b\geq\y\in\at{\gc\A}$. We want to prove $\x\rp\y\geq\c$,
so we also assume $\c\geq\z\in\at{\gc\A}$ and must now show $\x\rp\y\geq\z$.  Note that
$\x\leq\cov\x\leq\a$, $\y\leq\cov\y\leq\b$, and $\z\leq\cov\z\leq\c$. Also,
$\nog(\cov\x=\cov\y=\cov\z)$ since otherwise we would have $\a=\b=\c$, contradicting our
assumption.  Hence $[\cov\x,\cov\y,\cov\z]$ is not a 1-cycle, and is either a 2-cycle or
3-cycle.  Now $\eq$ contains all 2-cycles and 3-cycles by definition, so
$[\cov\x,\cov\y,\cov\z]$ is a cycle of $\eq$, hence $\cov\x\rp\cov\y\geq\cov\z$.  By
definition of splitting, $\x$ and $\y$ have a product equal to the product of their
covers in $\eq$, so we have $\x\rp\y=\cov\x\rp\cov\y$. If $\cov\x\neq\cov\y$ then
$\x\rp\y=\cov\x\rp\cov\y=\di\geq\c$ by Def.\,\ref{def1}. If $\cov\x=\cov\y$ then
$\a=\b$, hence $\a=\b\neq\c$ by the assumption $\nog(\a=\b=\c)$.  Therefore $\a\bp\c=0$
since $\a,\c\in\at{\gc\E}$. By Def.\,\ref{def1} we get
$\x\rp\y=\cov\x\rp\cov\y=\cov\x\rp\cov\x=\min{\cov\x}\geq\min{\a}\geq\c$.

To show part~(ii) of Def.\,\ref{spec}, we assume $\a\geq\x,\y\in\At{\gc\A}$ and
$\a\rp\a\geq\a$.  We wish to show that $\x\rp\y\bp\a\neq0$.  Note that $\a$ cannot be
an atom of $\eq$ because the assumption $\a\rp\a\geq\a$ fails for all atoms of $\eq$ by
Def.\,\ref{def1}.  Assume first that $\cov\x=\cov\y=\u$.  Since $\a$ is not an atom of
$\gc\E$, $\a$ is the join of two or more atoms of $\eq$, hence there is some atom
$\v\in\at\eq$ such that $\u\neq\v\leq\a$. We have $\x\rp\y=\cov\x\rp\cov\y=\u\rp\u$ by
the definition of splitting, but $\u\rp\u=\min\u\geq\v$ by Def.\,\ref{def1}, so
$0\neq\x\rp\y\bp\v\leq\x\rp\y\bp\a$, as desired.  Assume that $\cov\x\neq\cov\y$. Then
$\cov\x\rp\cov\y=\di$ by Def.\,\ref{def1}, so $\x\rp\y=\cov\x\rp\cov\y=\di\geq\a$, hence
$\x\rp\y\bp\a\neq0$.
\endproof
Lemma~\ref{lem1} suggests that we consider an arbitrary subalgebra $\gc\E$ of
$\eq$. Every subalgebra contains $\id$, but $\id$ is an atom in $\eq$, so it is an atom
in $\gc\E$ as well. Thus $\gc\E$ is integral, but $\gc\E$ is also symmetric since it is
the subalgebra of a symmetric algebra.  The diversity atoms of $\gc\E$ are disjoint and
join up to $\di$, so they partition the diversity atoms of $\eq$. In every relation
algebra, the relative product $\a\rp\b$ of \emph{distinct} diversity atoms $\a,\b$ of
$\gc\E$ is included in $\di$. On the other hand, in this case we have $\a\rp\b\geq\di$
because there are atoms $\x,\y$ of $\eq$ such that $\a\geq\x$, $\b\geq\y$, and
$\a\rp\b\geq\x\rp\y=\di$ by Def.\,\ref{def1}.  Every diversity atom of $\a\in\gc\E$
satisfies either $\a\rp\a=1$ or $\a\rp\a=\min\a$, for if $\a$ is an atom of $\eq$ (as
well as $\gc\E$) then $\a\rp\a=\min\a$ by Def.\,\ref{def1}, while if $\a$ is not an atom
of $\eq$, then it is the join of two or more atoms of $\eq$, say $\a\geq\cc1+\cc2$, so
\begin{align*}
\a\rp\a	&\geq(\cc1+\cc2)\rp(\cc1+\cc2)
\\	&=\cc1\rp\cc1+\cc1\rp\cc2+\cc2\rp\cc1+\cc2\rp\cc2
\\	&=\min{\cc1}+\di+\di+\min{\cc2}	&&\text{Def.\,\ref{def1}}
\\	&=1.
\end{align*}
Every subalgebra of $\eq$ can therefore be characterized by just two parameters:
$\alpha$, the number of diversity atoms $\a$ satisfying $\a\rp\a=\min\a$, and $\beta$,
the number of diversity atoms satisfying $\a\rp\a=1$.  The number of atoms in the
subalgebra is $1+\alpha+\beta$. The only restrictions on these parameters are
$\alpha+2\beta<\q$ and $0<\alpha+\beta$.  

An atom $\a$ of a symmetric integral relation algebra is said to be~\concept{flexible}
if $\a\rp\a=1$ and $\x\rp\a=\di$ for all diversity atoms $\x$ distinct from $\a$. Having
a flexible atom is a sufficient condition for representability; see
Comer~\cite[5.3]{Comer1984} or~\cite[Th.\,6]{Maddux1985}.  Since every proper subalgebra
of $\eq$ has at least one atom $\a$ satisfying $\a\rp\a=1$ and this atom is flexible by
Def.\,\ref{def1}, every proper subalgebra of $\eq$ is representable.

\section{\bf An infinite atom structure from two finite algebras}
In this section we use an arbitrary finite symmetric integral relation algebra $\gc\A$
and its subalgebra $\gc\E$ to construct a complete atomic algebra $\CEA\in\NA$ that has
subalgebras isomorphic to $\gc\A$ and $\gc\E$.

The~\concept{complex algebra} of the structure $\<\A,\C,\con\blank,\I\>$, where $\A$ is
a set, $\C$ is a ternary relation on $\A$, $\con\blank$ is a unary operation on $\A$,
and $\I\subseteq\A$, is the Boolean algebra of all subsets of $\A$ supplemented with
$\I$ as a distinguished element, the unary complex converse operation defined by
$\con\X=\{\con\x:\x\in\X\}$ for all $\X\subseteq\A$, and the complex relative
multiplication defined by
\begin{equation}\label{complexproduct}
	\X\rp\Y=\{\z:\exists\x\in\X,\,\exists\y\in\Y,\,\<\x,\y,\z\>\in\C\}
\end{equation}
for all $\X,\Y\subseteq\A$.  Every complete atomic $\NA$ is isomorphic to the complex
algebra of its atom structure~\cite[Th.3.13(2)]{MR662049}. Therefore, to define a
complete atomic $\NA$, as we do in the following definition, it is enough to describe
its atom structure.
\begin{definition}\label{def4}
  Assume $\gc\E\subseteq\gc\A\in\RA$ are finite symmetric integral relation
  algebras. Then $\CEA$ is the complete atomic $\NA$ with this atom structure: the atoms
  of $\CEA$ are $\id$ and the ordered pair $\pr\x\i$ for every diversity atom $\x$ of
  $\gc\A$ and every index $\i\in\omega$,
\begin{equation}\label{atoms.1}
	\at\CEA:=\{\id\}\cup\{\pr\x\i:\di\geq\x\in\at{\gc\A},\,\i\in\omega\},
\end{equation}
  the converse of every atom is itself, if $\T\subseteq\omega^3$ is defined for
  $\i,\j,\k\in\omega$ by
\begin{equation*}
	\Te\i\j\k\iff(\i\leq\j=\k)\,\lor\,(\j\leq\k=\i)\,\lor\,(\k\leq\i=\j),
\end{equation*}
  and $\cov\x$ is the atom of $\gc\E$ containing the atom $\x$ of $\gc\A$, then the
  cycles of $\CEA$ are, for all $\di\geq\x,\y,\z\in\at{\gc\A}$ and $\i,\j,\k\in\omega$,
\begin{align}
  &\label{1}	[\id,\id,\id],\ \ [\id,\pr\x\i,\pr\x\i],
\\&\label{3}	[\pr\x\i,\pr\y\j,\pr\z\k]
\quad\text{ if\,\, } \x\rp\y\geq\z\land\Big(\cov\x=\cov\y=\cov\z\ \Rightarrow\ \Te\i\j\k\Big).
\end{align}
  For any $\a\in\A$ and $\n\in\omega$, define the element $\an$ of $\CEA$ by
\begin{align}
  \an&=\sum\{\pr\x\i:\di\bp\a\geq\x\in\at{\gc\A},\,\n\leq\i\in\omega\}+\sum\{\id:\id\leq\a\},
\label{Jan}
\end{align}
  so if $\x$ is an atom of $\gc\A$, then
\begin{align}\label{atoman}
  \jn\x\n 
&=\begin{cases}	\sum\{\pr\x\i:\n\leq\i\in\omega\}	&\text{ if }\x\leq\di
	\\	\id					&\text{ if }\x=\id
  \end{cases},
\end{align}
\begin{enumerate}
\item\label{def4.iii}
  if $\di\geq\x\in\at{\gc\A}$, $\cov\x=\a\in\at{\gc\E}$, $\i,\j\in\omega$, and
  $\i\neq\j$, then
\begin{align*}
	\pr\x\i\rp\pr\x\i&=\jn{\di\bp\min\a\bp\x\rp\x}0+
	\join\{\pr\z\k:\k\leq\i,\,\a\bp\x\rp\x\geq\z\in\at{\gc\A}\}+\id,
\\	\pr\x\i\rp\pr\x\j&=\jn{\di\bp\min\a\bp\x\rp\x}0+
	\join\{\pr\z{\max(\i,\j)}:\a\bp\x\rp\x\geq\z\in\at{\gc\A}\},
\end{align*}
\item\label{def4.i}
  if $\di\geq\x,\y\in\at{\gc\A}$, $\i,\j\in\omega$, and $\cov\x\neq\cov\y$ then
\begin{align*}
	\pr\x\i\rp\pr\y\j
	&=\jn{\x\rp\y}0=\sum\big\{\pr\z\k:\x\rp\y\geq\z\in\at{\gc\A},\,\k\in\omega\big\},
\end{align*}
\item\label{def4.ii}
  if $\di\geq\x,\y\in\at{\gc\A}$, $\x\neq\y$, $\i,\j\in\omega$, $\i\neq\j$, and
  $\cov\x=\cov\y=\a\in\at{\gc\E}$ then
\begin{align*}
	\pr\x\i\rp\pr\y\j&=\jn{\di\bp\min\a\bp\x\rp\y}0+
	\join\{\pr\z{\max(\i,\j)}:\a\bp\x\rp\y\geq\z\in\at{\gc\A}\},
\\	\pr\x\i\rp\pr\y\i&=\jn{\di\bp\min\a\bp\x\rp\y}0+
	\join\{\pr\z\k:\k\leq\i,\,\a\bp\x\rp\y\geq\z\in\at{\gc\A}\}.
\end{align*}
\end{enumerate}
\end{definition}
Start with a finite symmetric integral relation algebra $\gc\A$ in which every atom is
splittable (in the sense of~\cite{MR1052567}). Let $\As\supseteq\gc\A$ be the relation
algebra obtained by splitting every atom $\a\in\at{\gc\A}$ into $\omega$ pieces
$\pr\a0,\pr\a1,\cdots$ so that $\a=\sum_{\i\in\omega}\pr\a\i$.  Splitting produces the
maximum set of cycles in the extension $\As\supseteq\gc\A$ that are consistent with
containing $\gc\A$ as a subalgebra. Let $\gc\E\subseteq\gc\A$ be a subalgebra of
$\gc\A$. From the atom structure of $\As$ we obtain a new atom structure whose complex
algebra is, in fact, isomorphic to $\CEA$, by deleting all the diversity cycles
$[\pr\a\i,\pr\b\j,\pr\c\k]$ of $\As$ which have the property that all the atoms in the
cycle lie below the same atom of $\gc\E$, and $\T(\i,\j,\k)$ fails to hold. This leaves
only a ``thin'' remnant of the cycles of $\As$ that we would classify as ``1-cycles of
$\gc\E$'' (because their atoms all lie below a single atom of $\gc\E$). The set of
1-cycles produced by splitting is significantly reduced by imposing the ``thinning
condition'' $\T(\i,\j,\k)$.  Those cycles of $\gc\A$ that are ``covered'' by 1-cycles of
$\gc\E$ are ``thinly reproduced'' in $\CEA$, while the 2- and 3-cycles of $\gc\A$ that
are covered by 2- or 3-cycles of $\gc\E$ are ``split'' into as many cycles as
possible. Treating 1-, 2-, and 3-cycles differently in various combinations, either
thinning or splitting each type of cycle, gives six more constructions that perhaps
should be examined with regard to Problem~\ref{prob}.
\begin{lemma}\label{lem2}
  $\gc\A$ is isomorphic, by $\a\mapsto\an$, to a subalgebra $\gc\A'$ of $\CEA$.
  \begin{equation*}\gc\A\cong\gc\A'\subseteq\CEA.\end{equation*}
\end{lemma}
\proof
Define the function $\varphi:\gc\A\to\CEA$ by $\varphi(\a)=\ao$ for all $\a\in\A$.  For
a key part of the proof that $\varphi$ embeds $\gc\A$ into $\CEA$, assume
$\di\geq\x,\y\in\at{\gc\A}$. We wish to prove that $\varphi(\x)\rp\varphi(\y)$ and
$\varphi(\x\rp\y)$ contain the same diversity atoms of $\CEA$.  (Proofs for the other
parts, involving preservation by $\varphi$ of the Boolean structure and identity
element, are fairly easy.)  

Consider an arbitrary diversity atom $\pr\z\k\in\at{\CEA}$, where
$\di\geq\z\in\at{\gc\A}$, $\k\in\omega$. Assume
$\pr\z\k\leq\varphi(\x)\rp\varphi(\y)$. Then there are $\u,\v\in\at{\CEA}$ such that
$\pr\z\k\leq\u\rp\v$, $\varphi(\x)\geq\u\in\at{\CEA}$, and
$\varphi(\y)\geq\v\in\at{\CEA}$.  By~\eqref{atoman} there are some $\i,\j\in\omega$ such
that $\u=\pr\x\i$ and $\v=\pr\y\j$. But then $[\pr\x\i,\pr\y\j,\pr\z\k]$ is a cycle of
$\CEA$, so $\x\rp\y\geq\z$ in $\gc\A$, which implies $\pr\z\k\leq\jn{\x\rp\y}0$, hence
$\pr\z\k\leq\varphi(\x\rp\y)$.  The argument is reversible.
\endproof
Every diversity atom in $\gc\A'$ is the join of an infinite set of atoms. Therefore
$\CEA$ cannot be a relation algebra if $\gc\A$ has any diversity atoms that are not
splittable. In fact, $\CEA$ satisfies all the axioms for relation algebras except
possibly the associative law, so $\CEA\in\NA$.  

Here is a computational lemma needed several times later.
\begin{lemma}\label{lemma}
  Assume $\gc\A$ is a special extension of $\gc\E$, $\a,\b$ are distinct diversity atoms
  of $\gc\E$, and $\u,\v$ are diversity atoms of $\CEA$.  If $\u\leq\jn\a0$ and
  $\v\leq\jn\b0$ then $\u\rp\v=\jn{\a\rp\b}0$. In particular, if $\a\rp\b=\di$, then
  $\u\rp\v=\di$.
\end{lemma}
\proof 
From $\u\leq\jn\a0$ and $\v\leq\jn\b0$ we have $\pr\x\i=\u$, $\pr\y\j=\v$,
$\x\leq\a=\cov\x$, and $\y\leq\b=\cov\y$, for some $\x,\y\in\at{\gc\A}$ and
$\i,\j\in\omega$.  The covers of $\x$ and $\y$ are different because
$\x\leq\cov\x=\a\neq\b=\cov\y\geq\y$. Rule Def.\,\ref{def4}\ref{def4.i} applies in this
case and says that $\pr\x\i\rp\pr\y\j=\jn{\x\rp\y}0$. Note that $\x\rp\y\leq\a\rp\b$.
Since $\gc\A$ is a special extension of $\gc\E$, we deduce from
Def.\,\ref{spec}\ref{spec.i} that every atom of $\gc\E$ below $\a\rp\b$ is also below
$\x\rp\y$, hence $\x\rp\y=\a\rp\b$.  We conclude that
$u\rp\v=\pr\x\i\rp\pr\y\j=\jn{\x\rp\y}0=\jn{\a\rp\b}0$.  If $\a\rp\b=\di$, then
$u\rp\v=\jn\di0$, but $\jn\di0$ is the diversity element $\di$ of $\CEA$, so
$\u\rp\v=\di$.
\endproof

\section{\bf Embedding Monk algebras}
Two elements of an atomic relation algebra are said to be~\concept{almost the same} if
their symmetric difference is the join of finitely many atoms.  We show in
Theorem~\ref{th1} below that if $\gc\A$ is a special extension of $\gc\E$ and $\gc\B$ is
the subalgebra of $\CEA$ generated by the atoms of $\CEA$, then the finitely generated
subalgebras of $\gc\B$ are finite and every element of $\gc\B$ almost the same as an
element of the subalgebra $\gc\E'$ of $\gc\B$ isomorphic to $\gc\E$ by
$\a\mapsto\jn\a0$.

For an example, suppose $\gc\A$ is a Monk algebra obtained from $\eq$ by splitting,
$4\leq\q\in\omega$, and $\gc\E$ is a subalgebra of $\eq$. By Lemma~\ref{lem1}, $\gc\A$
is a special extension of $\gc\E$, so Theorem~\ref{th1} applies to $\gc\A$ and $\gc\E$.
Next, we show in Theorem~\ref{th2}\ref{viii}\ref{ix} that if, in addition, $\gc\E$ has a
``flexible trio'' (Def~.\,\ref{def-trio} below) then $\gc\B$ is representable because
every finitely generated subalgebra of $\gc\B$ is included in a finite subalgebra of
$\gc\B$ that has the ``1-point extension property'' (Def~.\,\ref{def1pt} below).  In the
example, if $7\leq\q$ ($\gc\A$ has at least six colors) then $\eq$ has a subalgebra
$\gc\E$ with a flexible trio, so Theorem~\ref{th2}\ref{viii}\ref{ix} applies, and we
conclude that $\gc\B\in\RRA$.  Finally, we show in Theorem~\ref{th2}\ref{x} that if
$\gc\A$ has no 1-cycles then $\gc\B$ is not completely representable and $\CEA$, the
completion of $\gc\B$, is not representable.  Theorem~\ref{th2}\ref{x} applies to
$\gc\A$ because Monk algebras have no 1-cycles.  Cor.\,\ref{main} accordingly says that
every finite Monk algebra with six or more colors is a subalgebra of the
non-representable completion of an atomic representable relation algebra whose
finitely-generated subalgebras are finite.

The conclusion that $\CEA\notin\RRA$ can be obtained without Theorem~\ref{th2}\ref{x} in
case the Monk algebra $\gc\A$ is non-representable, which happens if the number of atoms
is large compared to the number of colors.  In this case the completion of $\gc\B$ is
non-representable simply because it has a non-representable subalgebra (isomorphic to
the non-representable Monk algebra $\gc\A$).

In Theorem~\ref{th1} we obtain conclusions just from knowing the extension
$\gc\E\subseteq\gc\A$ is special. Then in Theorem~\ref{th2} we also consider what
happens when, in addition, $\gc\E$ has a flexible trio and $\gc\A$ has no 1-cycles.
\begin{theorem}\label{th1}
  Assume $\gc\A$ and $\gc\E$ are finite symmetric integral relation algebras, $\gc\A$ is
  a special extension of $\gc\E\subseteq\gc\A$, and $\gc\B\subseteq\CEA$ is the
  subalgebra of $\CEA$ generated by $\at\CEA$.
\begin{enumerate}\item	\label{ii}
	$\gc\B$ is countable, atomic, symmetric, integral, and generated by its atoms.
\item	\label{ic}
	$\CEA$ and $\gc\B$ have the same atom structure.
\item	\label{ia}
	$\CEA$ is isomorphic to the complex algebra of the atom structure of $\gc\B$.
\item	\label{ib}
	$\CEA$ is the completion of $\gc\B$.
\item	\label{i}
	There are subalgebras $\gc\E'\subseteq\gc\A'\subseteq\CEA$ with
	$\gc\E'\cong\gc\E$ and $\gc\A'\cong\gc\A$.
\item	\label{vi}
	Every finitely generated subalgebra of $\gc\B$ is finite.
\item	\label{vii}
	Every element of $\gc\B$ is almost the same as an element of $\gc\E'$.
\end{enumerate}
\end{theorem}
\proof 
Parts~\ref{ii}--\ref{ib} require only the assumption that $\CEA$ is complete and atomic
and $\gc\B$ is the subalgebra of $\CEA$ generated by the atoms of $\CEA$.  Everything in
parts~\ref{ii}--\ref{ib} is either obvious or very easy to prove; see
\cite[Th.\,3.13]{MR662049} for part~\ref{ia}. Part~\ref{i} was proved in
Lemma~\ref{lem2}.  The assumption that $\gc\A$ is a special extension of $\gc\E$ is
needed only for Lemma~\ref{lem3} below, which is used to prove parts~\ref{vi}
and~\ref{vii}.
\begin{lemma}\label{lem3}
  For every $\n\in\omega$, $\bn$ is the set of atoms of a subalgebra of $\CEA$, where
  \begin{align}\label{Atn}
  \bn&=\{\id\}\cup\{\pr\x\i:\di\geq\x\in\at{\gc\A},\,\n>\i\in\omega\}
  \cup\{\an:\di\geq\a\in\at{\gc\E}\}.\end{align}
\end{lemma}
\proof
The elements of $\bn$ are disjoint and their join is $1$, so the set of joins of subsets
of $\bn$ is closed under the Boolean operations of $\CEA$ and, under those operations,
forms a Boolean algebra whose set of atoms is $\bn$.  The converse of everything in
$\bn$ is again in $\bn$ because conversion is the identity function on $\CEA$.  What
remains is to show the relative product $\u\rp\v$ of any two elements $\u,\v\in\bn$ is
the join of a subset of $\bn$.  For this it is enough to show that every element
$\w\in\bn$ is contained in or disjoint from $\u\rp\v$. This is clearly true whenever
$\u=\id$ or $\v=\id$ or $\w$ is itself an atom of $\CEA$, so we may assume $\w=\an$, for
some $\a\in\at{\gc\E}$, and $\u+\v\leq\di$.  We will show that if $\u\rp\v$ has nonempty
intersection with $\an$ then $\u\rp\v$ contains $\an$.

Suppose $\u\rp\v\bp\an\neq0$, where $\di\geq\u,\v\in\bn$, $\di\geq\a\in\at{\gc\E}$.
Then there are $\x,\y,\z\in\at{\gc\A}$ and $\i,\j,\k\in\omega$ such that
$\pr\x\i\leq\u$, $\pr\y\j\leq\v$, $\pr\z\k\leq\an$, and $\pr\x\i\rp\pr\y\j\geq\pr\z\k$.
For both cases below, note that $\cov\z=\a$ and $\n\leq\k$ by Def.\,\eqref{Jan}, and
$\x\rp\y\geq\z$ by~\eqref{3}.

{\bf Case~1.} $\nog(\cov\x=\cov\y=\cov\z)$.  From $\x\rp\y\geq\z$ we get
$\x\rp\y\bp\cov\z\neq0$ since $0\neq\z\leq\cov\z$, hence $\x\rp\y\geq\cov\z=\a$ by
Def.\,\ref{spec}\ref{spec.i}.  The implication in~\eqref{3} has a false hypothesis and
therefore holds trivially in this case for every atom of $\gc\A$ below $\a$. It follows
by~\eqref{3} and Def.\,\eqref{Jan} that $\pr\x\i\rp\pr\y\j\geq\ao\geq\an$.

{\bf Case~2.} $\cov\x=\cov\y=\cov\z=\a$.  In this case, by
$\pr\x\i\rp\pr\y\j\geq\pr\z\k$ and~\eqref{3} we have $\T(\i,\j,\k)$.  From
$\pr\x\i\leq\u\in\bn$ and the relevant definitions, it follows that if $\u$ is an atom
of $\CEA$ then $\u=\pr\x\i$ and $\i<\n$, while if $\u$ is not an atom of $\CEA$, then
$\pr\x\i\leq\u=\an$ since $\a=\cov\x$, and $\i\geq\n$.  Consequently, if both $\u$ and
$\v$ were atoms of $\CEA$, we would have $\i,\j<\n\leq\k$, contrary to
$\T(\i,\j,\k)$. Hence either $\u=\an$ or $\v=\an$. Since $\CEA$ is symmetric, these are
really the same case. We assume $\u=\pr\x\i$ and $\v=\an$, and will prove that
$\pr\x\i\rp\an\geq\an$.

Toward this end, assume $\w_\j\leq\an$ where $\n\leq\j\in\omega$ and $\w\leq\a=\cov\w$.
Now $\pr\x\i$ is in $\bn$, so $\i<\n\leq\j$, hence $\Te\i\j\j$.  From $\x\rp\y\geq\z$
and $\cov\x=\cov\y=\cov\z=\a$ we get $\a\rp\a\bp\a\neq0$, but $\x\leq\a$ and $\w\leq\a$,
so $\x\rp\w\bp\a\neq0$ by Def.\,\ref{spec}\ref{spec.ii}.  We may therefore choose an atom
$\t\in\gc\A$ such that $\t\leq\x\rp\w\bp\a$.  Then $\cov\t=\a$ and $\Te\i\j\j$, so
$[\pr\x\i,\pr\t\j,\pr\w\j]$ is a cycle of $\CEA$ by~\eqref{3}, and $\pr\t\j\leq\an$
since $\t\leq\a$ and $\n\leq\j$, so $\w_\j\leq\pr\x\i\rp\pr\t\j\leq\pr\x\i\rp\an$.
Since this holds for all atoms $\w_\j$ below $\an$, we have proved
$\an\leq\pr\x\i\rp\an$.

We have shown that every product of two elements of $\bn$ is the join of a subset of
$\bn$.  It follows that $\u\rp\v=\tsum\{\w:\u\rp\v\geq\w\in\bn\}$ for all
$\u,\v\in\bn$. Hence,for all $\U,\V\subseteq\bn$, we have
\begin{align*}
\tsum\U\rp\tsum\V
	&=\tsum\{\u\rp\v:\u\in\U,\,\v\in\V\}
\\	&=\tsum\{\tsum\{\w:\u\rp\v\geq\w\in\bn\}:\u\in\U,\,\v\in\V\}
\\	&=\tsum\{\w:\u\rp\v\geq\w\in\bn,\,\u\in\U,\,\v\in\V\}
\\	&\in\{\tsum\X:\X\subseteq\bn\}.
\end{align*}
Therefore $\{\tsum\X:\X\subseteq\bn\}$ is closed under relative multiplication and is a
subalgebra of $\CEA$.
\endproof
We return to the proof of Theorem~\ref{th1}. Suppose $\gc\F$ is a finitely generated
subalgebra of $\gc\B$. Since $\gc\B$ is itself generated by $\at{\CEA}$, there is a
finite set of atoms $\X\subseteq\at{\CEA}$ such that $\gc\F$ is contained in the
subalgebra of $\CEA$ generated by $\X$. Since $\X$ is finite and
$\at{\CEA}\subseteq\bigcup_{\n\in\omega}\bn$, we may choose a sufficiently large
$\n\in\omega$ so that $\X\subseteq\bn$. Then $\gc\F$ is contained in the subalgebra
$\gc\B_\n$ of $\CEA$ generated by $\bn$. The subalgebra $\gc\B_\n$ is finite by the
lemma, since its set of atoms is the finite set $\bn$, so $\gc\F$ is also
finite. Hence~\ref{vi} holds.  This argument also shows that every element of $\gc\B$
is, for some $\n\in\omega$, included in a subalgebra whose set of atoms is $\bn$, and
hence is a join of elements of $\bn$. But every join of elements of $\bn$ is almost the
same as one of the atoms of $\gc\E'$.  Hence~\ref{vii} holds.
\endproof
As it happens, every finitely-generated subalgebra of $\CEA$ (not just $\gc\B$) is also
finite, even if the extension $\gc\E\subseteq\gc\A$ is not special. To prove this, one
argues that for every finite subset $\F$ of $\CEA$ there is some $\n\in\omega$ and some
finite partition $\mathcal{P}$ of $\{\i:\n\leq\i\in\omega\}$ such that
\begin{equation*}
	\{\id\}
	\cup\{\pr\x\i:\x\in\at{\gc\A},\,\n>\i\in\omega\}
	\cup\{\textstyle\sum\{\pr\x\i:\i\in\P\}:\x\in\at{\gc\A},\,\P\in\mathcal{P}\}
\end{equation*}
is the set of atoms of a subalgebra of $\CEA$ that contains $\F$. The remaining details
of this proof are omitted since this fact is not needed and it is also not in itself
enough to prove Lemma~\ref{lem3}. On the other hand, there is a special case which is
easy to prove and needed later.
\begin{lemma}\label{lem3a}
  Assume $\gc\A\supseteq\gc\E$ are finite symmetric integral relation algebras. For
  every $\n\in\omega$,
\begin{align}\label{Atn.1}
&\{\id\}\cup\{\pr\x\i:\di\geq\x\in\at{\gc\A},\,\n>\i\in\omega\}
\cup\{\an:\di\geq\a\in\at{\gc\A}\}.
\end{align}
is the set of atoms of a subalgebra of $\CEA$.
\end{lemma}
\proof	The proof is similar to, but simpler, than the proof of Lemma~\ref{lem3}. The
	closure of the set of joins of subsets of \eqref{Atn.1} under relative
	multiplication is an immediate consequence of
	Def.\,\ref{def4}\ref{def4.iii}\ref{def4.i}\ref{def4.ii}.
\endproof
For the next theorem we need some definitions. A relation algebra has the 1-point
extension property if, loosing speakly, every ``finite partial representation'' $\mo$
can be extended by one point wherever this is needed. We make this precise as follows.
\begin{definition}\label{def1pt}
  For any $\k\in\omega$ and any atomic relation algebra $\gc\A$, $\B_\k(\gc\A)$ is the
  set of functions $\mo:\k\times\k\to\at{\gc\A}$ that satisfy the following conditions.
\begin{enumerate}
\item[{\rm (B$_0$)}] $\mo_{\i,\i}\leq\id$ for all $\i<\k$,
\item[{\rm (B$_1$)}] $\con\mo_{\i,\j}=\mo_{\j,\i}$ for all $\i,\j<\k$,
\item[{\rm (B$_2$)}] $\mo_{\i,\l}\rp\mo_{\l,\j}\geq\mo_{\i,\j}$ for all $\i,\j,\l<\k$.
\end{enumerate}
The elements of $\B_\k(\gc\A)$ are called \concept{basic matrices}.  A matrix $\mo$
satisfies the \concept{identity condition} if $\mo_{\l,\m}=\id$ iff $\l=\m$ for all
$\l,\m<\k$.  We say that $\gc\A$ has the \concept{1-point extension property} if,
assuming $\mo\in\B_\k(\gc\A)$, $\mo$ satisfies the identity condition, $\x,\y$ are
diversity atoms of $\gc\A$, $\i,\j<\k$, $\i\neq\j$, and $\mo_{\i,\j}\leq\x\rp\y$, there
are basic matrix $\mo'\in\B_{\k+1}(\gc\A)$ satisfying the identity condition such that
$\mo_{\l,\m}=\mo'_{\l,\m}$ for all $\l,\m<\k$, $\mo'_{\i,\k}=\x$, and $\mo'_{\k,\j}=\y$.
\end{definition}
\begin{definition}\label{def-cr}
A relation algebra $\gc\A$ is \concept{completely representable} if it has a complete
representation, where a representation $\rho$, mapping $\gc\A$ into some algebra of
binary relations, is \concept{complete} if it preserves all joins, \ie, if $\X$ is a
subset of $\gc\A$ whose join $\join\X$ exists in $\gc\A$, then
$\rho\(\join\X\)=\bigcup_{\x\in\X}\rho(\x)$.
\end{definition}
\begin{definition}\label{def-trio}
  Three diversity atoms $\a,\b,\c$ of a symmetric integral relation algebra $\gc\A$ are
  said to be a \concept{flexible trio} if
\begin{align}
\label{m2}	\a\rp\a	&=\b\rp\b=\c\rp\c=1,
\\
\label{m1}	\a\rp\b	&=\a\rp\c=\b\rp\c=\di,
\intertext{and, for every atom $\x\notin\{\id,\a,\b,\c\}$,}
\label{m}	\x\rp\a	&=\x\rp\b=\di\,\,\lor\,\,
		\x\rp\a=\x\rp\c=\di\,\,\lor\,\,\x\rp\b=\x\rp\c=\di.
\end{align}
\end{definition}
Theorem~\ref{th3} (relegated to an Appendix) shows that having a flexible trio is
sufficient for representability.

If each of $\a$, $\b$, and $\c$ is a flexible atom then $\a,\b,\c$ is a flexible
trio. Hence any three diversity atoms of $\eq$ form a flexible trio. It can happen that
$\a,\b,\c$ is a flexible trio but none of $\a,\b,\c$ is flexible, as in, for example,
the symmetric integral relation algebra with seven atoms $\id,\a,\b,\c,\d,\e,\f$ and all
diversity cycles \emph{except} $[\a,\d,\d]$, $[\b,\e,\e]$, and $[\c,\f,\f]$. This
algebra has no flexible atoms.

The following theorem applies to all Monk algebras with at least six colors, but many
other algebras also satisfy its hypotheses. For an example, let $\gc\A$ be the symmetric
integral relation algebra whose atoms are $\id$, $\a_1$, $\a_2$, $\a_3$, $\b_1$, $\b_2$,
$\b_3$, $\c_1$, $\c_2$, and $\c_3$, and whose diversity cycles consist of none of the
1-cycles, all of the 2-cycles, and all the 3-cycles \emph{except} $[\a_1,\a_2,\a_3]$,
$[\b_1,\b_2,\b_3]$, and $[\c_1,\c_2,\c_3]$.  Let $\gc\E$ be the subalgebra of $\gc\A$
whose atoms are $\id$, $\a=\a_1+\a_2+\a_3$, $\b=\b_1+\b_2+\b_3$, and
$\c=\c_1+\c_2+\c_3$.  Then $\gc\A$ is a special extension of $\gc\E$, $\a,\b,\c$ is a
flexible trio of $\gc\E$, and $\gc\A$ has no 1-cycles, so Th.~\ref{th2} can applied to
conclude that $\gc\A$ is in the class $\K$ of Problem~\ref{prob}, but $\gc\A$ is not a
Monk algebra.  More than 3000 additional examples can be obtained by deleting any or all
of the following 2-cycles: $[\a_1,\a_2,\a_2]$, $[\a_2,\a_3,\a_3]$, $[\a_1,\a_1,\a_3]$,
$[\b_1,\b_2,\b_2]$, $[\b_2,\b_3,\b_3]$, $[\b_1,\b_1,\b_3]$, $[\c_1,\c_2,\c_2]$,
$[\c_2,\c_3,\c_3]$, $[\c_1,\c_1,\c_3]$, and restoring any or all of the deleted
3-cycles, but at least one of the 2- or 3-cycles must be deleted to insure $\gc\A$ is
not a Monk algebra. This scheme retains enough 2-cycles that the extension is always
special.
\begin{theorem}\label{th2}
  Assume $\gc\E\subseteq\gc\A$ are finite symmetric integral relation algebras, $\gc\A$
  is a special extension of $\gc\E$, and $\gc\E$ has a flexible trio.  If $\gc\B$ is the
  subalgebra of $\CEA$ generated by $\at\CEA$, then
\begin{enumerate}
\item	\label{viii}
	every finitely generated subalgebra of $\gc\B$ is contained in a subalgebra of
	$\gc\B$ that has the 1-point extension property,
\item	\label{ix}
	$\gc\B$ is representable,
\item	\label{x}
	if $\gc\A$ has no 1-cycles, \ie, $\u\rp\u\bp\u=0$ whenever
	$\id\neq\u\in\at{\gc\A}$, then $\gc\B$ is not completely representable and the
	completion of $\gc\B$ is not representable.
\end{enumerate}
\end{theorem}
\proof
Suppose $\gc\F$ is a finitely generated subalgebra of $\gc\B$.  By the argument at the
end of the proof of Th.\,\ref{th1}, there is some $\n\in\omega$ such that $\gc\F$ is
contained in the subalgebra $\gc\B_\n\subseteq\gc\B$ with atoms $\at{\gc\B_\n}=\bn$.
Let $\a,\b,\c$ be a flexible trio of $\gc\E$. We will show that
$\jn\a\n,\jn\b\n,\jn\c\n$ is a flexible trio of $\gc\B_\n$. Consider the product of the
first two elements of the trio.  Note that $\jn\a\n\rp\jn\b\n\leq\di$ since $\jn\a\n$
and $\jn\b\n$ are disjoint atoms of $\gc\B_\n$.  We have
\begin{align*}
  & \jn\a\n\rp\jn\b\n
\\&=\sum\{\u\rp\v:\jn\a\n\geq\u\in\at\CEA,\,\jn\b\n\geq\v\in\at\CEA\}
\end{align*}
but every disjunct $\u\rp\v$ in this last join is $\di$ by Lemma~\ref{lemma} and the
assumption $\a\rp\b=\di$, so $\jn\a\n\rp\jn\b\n=\di$.  Similarly,
$\jn\a\n\rp\jn\c\n=\di=\jn\b\n\rp\jn\c\n$. Thus~\eqref{m1} holds.

For \eqref{m}, consider a diversity atom of $\gc\B_\n$ that is not one of
$\jn\a\n,\jn\b\n,\jn\c\n$. It is either an atom of $\CEA$ or at atom of $\gc\B_\n$ with
the form $\jn\d\n$, where $\d$ is a diversity atom of $\gc\E$ distinct from $\a,\b,\c$.

We first consider $\jn\d\n$.  Now $\d$ multiplies to $\di$ with two of $\a,\b,\c$
by~\eqref{m}, say $\a\rp\d=\b\rp\d=\di$. Choose atoms $\x,\y$ of $\gc\A$ with $\x\leq\a$
and $\y\leq\d$. Then $\jn\a\n\rp\jn\d\n\leq\di$ since $\jn\a\n$ and $\jn\d\n$ are
disjoint, and $\jn\a\n\rp\jn\d\n\geq\pr\x\n\rp\pr\y\n$, but $\pr\x\n\rp\pr\y\n=\di$ by
Lemma~\ref{lemma} because $\a\rp\d=\di$, so $\jn\a\n\rp\jn\d\n=\di$.  Similarly
$\jn\b\n\rp\jn\d\n\geq\di$, so the atom $\jn\d\n$ multiplies to $\di$ with two of
$\jn\a\n,\jn\b\n,\jn\c\n$, as desired.

Next consider at atom $\u$ of $\CEA$. It has the form $\u=\pr\x\i$ for some diversity
atom $\x$ of $\gc\A$ and some $\i<\n$.  We claim that the product of $\cov\x$ with (at
least) two elements in the trio $\a,\b,\c$ is $\di$, say $\a\rp\cov\x=\b\rp\cov\x=\di$.
This follows from~\eqref{m} if $\cov\x$ is a diversity atom distinct from $\a,\b,\c$,
but if $\cov\x$ is one of $\a,\b,\c$, then it follows from~\eqref{m1}.  Choose an atom
$\a\geq\y\in\at{\gc\A}$. Then $\pr\x\i\rp\jn\a\n\geq\pr\x\i\rp\pr\y\n=\di$ by
$\a\rp\d=\di$ and Lemma~\ref{lemma}.  Similarly $\pr\x\i\rp\jn\b\n=\di$, so the atom
$\u=\pr\x\i$ multiplies to $\di$ with two of $\jn\a\n,\jn\b\n,\jn\c\n$, as desired. This
finishes the proof of~\eqref{m} for $\jn\a\n,\jn\b\n,\jn\c\n$.

For~\eqref{m2}, we will prove $\jn\a\n\rp\jn\a\n=1$ from $\a\rp\a=1$.  Assume
$\u=\pr\x\i\in\at\CEA$, $\di\geq\x\in\at{\gc\A}$, and $0\leq\i\in\omega$.  Then
$\x\leq1=\a\rp\a=\sum\{\y\rp\z:\a\geq\y,\z\in\at{\gc\A}\}$ so there are atoms
$\a\geq\y,\z\in\at{\gc\A}$ such that $\x\leq\y\rp\z$.  Note that $\cov\y=\cov\z=\a$.
Choose any $\j$ such that $\max(\i,\n)\leq\j\in\omega$.  Then $\T(\i,\j,\j)$ holds, so
$[\pr\x\i,\pr\y\j,\pr\z\j]$ is a cycle of $\CEA$ by~\eqref{3}, hence
$\u=\pr\x\i\leq\pr\y\j\rp\pr\z\j\leq\jn\a\n\rp\jn\b\n$. This shows
$\jn\a\n\rp\jn\a\n=1$, and we obtain $\jn\b\n\rp\jn\b\n=1=\jn\c\n\rp\jn\c\n$ similarly
from $\b\rp\b=\c\rp\c=1$.

This completes the proof that $\jn\a\n,\jn\b\n,\jn\c\n$ is a flexible trio of
$\gc\B_\n$. By Theorem~\ref{th3} below, $\gc\B_n$ has the 1-point extension property and
is therefore representable. Every finitely generated subalgebra of $\gc\B$ is
representable, hence $\gc\B$ is representable since $\RRA$ is a variety.  Thus
parts~\ref{viii} and~\ref{ix} hold.

For part \ref{x}, assume such that $\u\rp\u\bp\u=0$ whenever $\id\neq\u\in\at{\gc\A}$.
Suppose that $\rho$ is a complete representation of $\gc\B$. By definition, a complete
representation preserves all joins. In particular, the join of the diversity atoms is
$\di$, so the union of the representations of the diveristy atoms must be the diversity
relation, consisting of all pairs of distinct objects. Thus every pair $\<\i,\j\>$ with
$\i\neq\j$ is in the representation of a diversity atom of $\gc\B$. But
$\gc\A\cong\gc\A'\subseteq\gc\B$ via $\a\mapsto\jn\a0$, so the representation $\rho$
determines a coloring of the edges of $\K_\omega$ (the complete graph on countably many
vertices) as follows: the ``color'' of the edge $(\i,\j)$ is $\a\in\at{\gc\A}$ if
$\<\i,\j\>\in\rho(\jn\a0)$.  This coloring has no monochrome triangles because $\gc\A$
has no $1$-cycles, contrary to Berge~\cite[Prop.\,1, p.\,440]{MR0357172}.  Therefore
there is no such $\rho$ and $\gc\B$ is not completely representable.
\endproof
\begin{cor}\label{main}
If $\gc\A$ is a finite Monk algebra with six or more colors then $\gc\A$ is in the class
$\K$ defined in Problem~1. In fact, $\gc\A$ is a subalgebra of the completion of a
relation algebra $\gc\B$ such that
\begin{itemize}
\item	$\gc\B$ is a countable, atomic, symmetric, integral relation algebra that is
	generated by its atoms,
\item	every finitely generated subalgebra of $\gc\B$ is contained in a
	finite subalgebra of $\gc\B$ with the 1-point extension property,
\item	the completion of $\gc\B$ has the same atom structure as $\gc\B$, is isomorphic
	to the complex algebra of the atom structure of $\gc\B$, and is not
	representable,
\item	$\gc\B$ is representable but not completely representable.
\end{itemize}
\end{cor}
The smallest example to which these considerations apply is $\gc\A=\ex$, the example
considered earlier.  This algebra is a Monk algebra with six colors and no 1-cycles,
obtained from itself by splitting.  By the Street-Whitehead-Comer method, $\ex\in\RRA$
because $\ex$ has square representations on sets containing 97, 157, and 277
elements. [AMS Meeting, Iowa City, March 2011.] Suppose the diversity atoms of $\ex$ are
$\cc1$, \dots, $\cc6$. Let $\aa1=\cc1+\cc2$, $\aa2=\cc3+\cc4$, and
$\aa3=\cc5+\cc6$. Then $\gc\A$ is special extension of the subalgebra $\gc\E$ whose
atoms are $\id,\aa1,\aa2,\aa3$, and $\aa1,\aa2,\aa3$ is a flexible trio of individually
flexible atoms in $\gc\E$.  Then $\RRA$ is not closed under completions because $\Cex$
is the non-representable completion of the atomic representable subalgebra of $\Cex$
generated by $\at\Cex$. In the next section we compute the exact degree of
non-representability of $\Cex$.

\section{\bf Cylindric algebras}
$\CA\n$ is the class of $\n$-dimensional cylindric algebras.  Given a cylindric algebra
$\gc\D\in\CA\n$ of dimension $\n\geq3$, the~\concept{relation algebraic reduct}
$\Ra{\gc\D}$ is defined in~\cite[Def.\,5.3.7]{MR781930} and is a relation algebra if
$\n\geq4$ by~\cite[Def.\,5.3.8]{MR781930}.  For any class $\K\subseteq\CA\n$ with
$3\leq\n$, let $\mathsf{Ra}\K$ be the class of relation-algebraic reducts of subalgebras
of neat 3-dimensional reducts of algebras in $\K$:
\begin{equation*}
	\mathsf{Ra}\K=\gc{Ra}^*\Sub\Nrr3\K
\end{equation*}
By~\text{\cite[5.3.9, 5.3.16, 5.3.17]{MR781930}}, we have
\begin{align*}
  \RRA&=\prod_{\n\in\omega}\mathsf{Ra}\CA{\n+4}
	\subseteq\cdots\subseteq
	\mathsf{Ra}\CA5\subseteq\mathsf{Ra}\CA4=\RA.
\end{align*}
Every non-representable relation algebra lies somewhere on this chain.  The location of
the example $\Cex$ is determined by the main result in this section, which implies
\begin{equation}\label{cut}
	\Cex\in\mathsf{Ra}\CA7\minus\mathsf{Ra}\CA8.
\end{equation}
\begin{definition}\label{def-basis}
Assume $\gc\A\in\NA$ is atomic and $\k\leq\omega$.  Two basic matrices $\mo$ and $\mo'$
in $\B_\k(\gc\A)$ \concept{agree up to} $\i$ if $\mo_{\l,\m}=\mo'_{\l,\m}$ whenever
$\i\neq\l,\m\in\k$, and they \concept{agree up to} $\i,\j$ if $\mo_{\l,\m}=\mo'_{\l,\m}$
whenever $\i,\j\neq\l,\m\in\k$. We say that $\MM\subseteq\B_\k(\gc\A)$ is an
$\k$-\concept{dimensional relational basis for} $\gc\A$ if
\begin{enumerate}
\item[{\rm(R$_0$)}]
	for every atom $\a\in\at{\gc\A}$ there is a basic matrix $\mo\in\MM$ such that
	$\mo_{0,1}=\a$,
\item[{\rm(R$_1$)}]
	if $\mo\in\MM$, $\i,\j<\k$, $\x,\y\in At\gc\A$, $\mo_{\i,\j}\leq\x\rp\y$, and
	$\i,\j\neq\l<\k$, then there is some $\mo'\in\MM$ such that $\mo$ and $\mo'$
	agree up to $\l$, $\mo'_{\i,\l}=\x$, and $\mo'_{\l,\j}=\y$.
\end{enumerate}
For any $\i,\j<\k$ let
\begin{align*}
	\T^\k_\i(\gc\A)&=\lb\<\mo,\mo'\>\in\B_\k(\gc\A)\times\B_\k(\gc\A): 
	\text{$\mo$ and $\mo'$ agree up to $\i$}\rb,
\\	\E^\k_{\i,\j}(\gc\A)&=\lb\mo\in\B_\k(\gc\A):\mo_{\i,\j}\leq\id\rb.
\end{align*}
We say that $\MM\subseteq\B_\k\gc\A$ is a \concept{$\k$-\nobreak dimensional cylindric
basis for} $\gc\A$ if
\begin{enumerate}
\item[{\rm(C$_0$)}]
if $\a,\b,\c\in\at{\gc\A}$, and $\a\leq\b\rp\c$, then there is a basic matrix $\mo\in\MM$
such that $\mo_{01}=\a$, $\mo_{02}=\b$, and $\mo_{21}=\c$,
\item[{\rm(C$_1$)}]
if $\mo,\mo'\in\MM$,\ $\i,\j<\k$,\ $\i\neq\j$, and $\mo$ agrees with $\mo'$ up to
$\i,\j$, then there is some $\mo''\in\MM$ such that $\mo''$ agrees with $\mo$ up to
$\i$, and $\mo''$ agrees with $\mo'$ up to $\j$, \ie, $\<\mo'',\mo\>\in\T^\k_\i(\gc\A)$
and $\<\mo'',\mo'\>\in\T^\k_\j(\gc\A)$,
\item[{\rm(C$_2$)}]
if $\mo\in\MM$ and $\i,\j<\k$ then $\mo[\i/\j]\in\MM$, where $[\i/\j](\m)=\m$ if
$\i\neq\m<\k$, and $[\i/\j](\i)=\j$.
\end{enumerate}
For every $\MM\subseteq\B_\k(\gc\A)$, let
\begin{equation}\label{complexalg}
	\Ca\MM=\Cm{\<\MM,\T_\i,\E_{\i\j}\>_{\i,\j<\k}}
\end{equation}
be the complex algebra of the relational structure $\<\MM,\T_\i,\E_{\i\j}\>_{\i,\j<\k}$,
as defined in~\cite[2.7.33]{MR0314620}, where $\E_{\i\j}=\E^\k_{\i\j}(\gc\A)\cap\MM$ and
$\T_\i=\T^\k_\i(\gc\A)\cap\(\MM\times\MM\)$ for all $\i,\j<\k$.
\end{definition}
\begin{theorem}\label{th4}
Assume $4\leq\r\in\omega$, $\gc\A=\ep$, and the atoms of $\gc\A$ are
\begin{equation*}	\cc0=\id,\,\cc1,\,\cc2,\,\cc3,\,\cc4,
	\,\cc5,\,\cc6,\,\cc7,\,\cdots,\,\cc{\r+2}.
\end{equation*}
Then $\gc\A$ is a special extension of a subalgebra $\gc\E$ whose $\r$ atoms are
\begin{align*}
   \aa0&=\cc0=\id&\aa1&=\cc1+\cc2&\aa2 &=\cc3+\cc4&\aa3     &=\cc5+\cc6
\\     &         &\aa4&=\cc7     &     &\cdots    &\aa{\r-1}&=\cc{\r+2}
\end{align*}
and $\aa1,\aa2,\aa3$ is a flexible trio of $\gc\E$, so by Ths.1,2, the atom-generated
subalgebra of the complete atomic relation algebra $\CEA$ is an atomic atom-generated
symmetric integral representable relation algebra with finite finitely-generated
subalgebras, and if $3\leq\n\leq\r+3$ then
\begin{enumerate}
\item	$\B_\n(\CEA)$ is an $\n$-dimensional cylindric basis for $\CEA$,
\label{li}
\item	$\Ca{\B_\n(\CEA)}$ is a complete atomic $\n$-dimensional cylindric algebra,
\label{lii} 
\item	$\CEA$ is isomorphic to the relation algebraic reduct of $\Ca{\B_\n(\CEA)}$ and
\label{liii} $\CEA\in\mathsf{Ra}{\CA\n}$,
\item	$\Ca{\B_\n(\CEA)}\notin\Sub\nr_\n\CA{\r+4}$.
\label{liv}
\end{enumerate}
\end{theorem}
\proof
By~\cite[Th.\,7]{MR1011183}, in order to prove $\B_\n(\CEA)$ is a cylindric basis for
$\CEA$ it is enough to show, given $\n-2$ pairs of diversity atoms
$\uu1,\vv1,\cdots,\uu{\n-2},\vv{\n-2}$ of $\CEA$, that
$\prod_{1\leq\i\leq\n-2}\uu\i\rp\vv\i\neq0$.  We will find a diversity atom $\w$, such
that $\w$ is included in every product $\uu\i\rp\vv\i$, $1\leq\i\leq\n-2$.  Any product
$\uu\i\rp\vv\i$ that is equal to $\di$ or $1$ imposes no restriction on our choice of
$\w$. We therefore assume that~\emph{none} of the products is $\di$ or $1$, \ie,
$\di\neq\uu\i\rp\vv\i\neq1$ whenever $1\leq\i\leq\n-2$.

Consequently, for every product $\uu\i\rp\vv\i$ we know that there cannot
be~\emph{distinct} atoms $\a,\b\in\at{\gc\E}$ such that $\uu\i\leq\jn\a0$ and
$\vv\i\leq\jn\b0$, because we would obtain $\uu\i\rp\vv\i=\jn{\a\rp\b}0$ from $\a\neq\b$
by Lemma~\ref{lemma}, and a computation in $\gc\E$ shows $\a\rp\b=\di$ since $\a\neq\b$,
forcing $\uu\i\rp\vv\i=\di$ in $\CEA$, contrary to our assumption that no product is
$\di$ or $1$.  Therefore, there is a function $\f:\{1,\cdots,\n-2\}\to\{1,\cdots,\r-1\}$
such that
\begin{equation}\label{*}
\uu\i+\vv\i\leq\jn{\aa{\f(\i)}}0\text{ for all }\i\in\{1,\cdots,\n-2\}.
\end{equation}
We see next that every index in $\{1,\cdots,\r-1\}$ is in the range of $\f$.  Suppose
some index $\j\{1,\cdots,\r-1\}$ is not in the range of $\f$.  Consider any product
$\uu\i\rp\vv\i$ with $1\leq\i\leq\n-2$.  Let $\k=\f(\i)$ and note that $\k\neq\j$.
There are atoms $\x,\y\in\at{\gc\A}$ such that
\begin{align}\label{8a}
	\x+\y\leq\aa\k=\cov\x=\cov\y,
&&	\uu\i\leq\jn\x0,
&&	\vv\i\leq\jn\y0.
\end{align}
Since $\k\neq\j$, $\min{\aa\k}\geq\aa\j$.  By Def.\,\ref{def4}\ref{def4.iii}\ref{def4.ii}
and~\eqref{8a}, $\uu\i\rp\vv\i\geq\jn{\di\bp\min{\aa\k}\bp\x\rp\y}0$.  If $\x\neq\y$ then
$\x\rp\y=\di$ in $\gc\A$, so
\begin{equation}\label{have}
	\uu\i\rp\vv\i\geq\jn{\di\bp\min{\aa\k}\bp\x\rp\y}0
	=\jn{\di\bp\min{\aa\k}}0\geq\jn{\aa\j}0.
\end{equation}
If $\x=\y$ then $\x\rp\y=\min\x\geq\min{\aa\k}$ in $\gc\A$, so
$\min{\aa\k}\bp\x\rp\y=\min{\aa\k}$, and again we have~\eqref{have}.\footnote{By the
way, we've shown $\j\neq\f(\i)\implies\uu\i\rp\vv\i\geq\jn{\aa\j}0$.}
Since~\eqref{have} holds for every $\i$, we obtain much more than
$\prod_{1=\i}^{\n-2}\uu\i\rp\vv\i\neq0$, in fact,
\begin{equation}
	0\neq\jn{\aa\j}0\leq\prod_{\i=1}^{\n-2}\uu\i\rp\vv\i.
\end{equation}
Therefore, assume that $\f$ is surjective. Next we show that $\f$ is actually maps two
distinct indices onto each of the atoms $\aa1$, $\aa2$, and $\aa3$, \ie, those atoms of
$\gc\E$ that are the join of two atoms of $\gc\A$.  We prove this only for $\aa1$.
Since $1$ is in the range of $\f$, we'll suppose, for specificity and simplicity of
notation, that $1=\f(1)$, \ie, $\uu1+\vv1\leq\jn{\aa1}0$. We wish to show that
$1=\f(\i)$ for some $\i\neq1$, so we assume this does \emph{not} happen, \ie, assume
$1\neq\f(\i)$ for all $\i\in\{2,\cdots,\n-2\}$.  Now $\aa1=\cc1+\cc2$, so there are
atoms $\x,\y\in\{\cc1,\cc2\}$ and indices $\k,\l\in\omega$ such that $\x+\y\leq\aa1$,
$\uu1=\pr\x\k$, and $\vv1=\pr\y\l$. Let $\m=\max(\k,\l)$. Notice that $\T(\k,\l,\m)$
holds and $\ep$ contains the 2-cycles $[\cc1,\cc1,\cc2]$ and $[\cc1,\cc2,\cc2]$.  It
follows by~\eqref{3} that
\begin{align*}
  	\pra1\k\rp\pra1\l			&\geq\pra2\m
\\	\pra1\k\rp\pra2\l=\pra2\k\rp\pra1\l	&\geq\pra1\m+\pra2\m
\\	\pra2\k\rp\pra2\l			&\geq\pra1\m
\end{align*}
We may therefore let $\w=\pra1\m$ if $\x=\y=\cc1$, $\w=\pra2\m$ if $\x=\y=\cc2$, and
either $\w=\cc1$ or $\w=\cc2$ if $\x\neq\y$. In every case, $\uu1\rp\vv1\geq\w$.  For
products other than $\uu1\rp\vv1$, note that if $2\leq\i\leq\n-2$, then by our
assumption we have $1\neq\f(\i)$, hence by the footnote,
$\uu\i\rp\vv\i\geq\jn{\aa1}0\geq\w$. This shows
$\w\leq\prod_{\i=2}^{\n-2}\uu\i\rp\vv\i$, which, together with $\w\leq\uu1\rp\vv1$,
gives us $\w\leq\prod_{\i=1}^{\n-2}\uu\i\rp\vv\i$.

At this point we know that either we are done because we have proved
$\prod_{\i=1}^{\n-2}\uu\i\rp\vv\i\neq0$, or else $\f$ maps at least one index from
$\{1,\cdots,\n-2\}$ onto each of the indices in $\{4,\cdots,\r-1\}$, and $\f$ maps at
least two indices from $\{1,\cdots,\n-2\}$ onto each of the indices in $\{1,2,3\}$, But
$|\{1,\cdots,\n-2\}|=\n-2$, $|\{4,\cdots,\r-1\}|=\r-4$, and $|\{1,2,3\}|=3$, so we must
have $\n-2\geq\r-4+2\cdot3=\r+2$, but our restriction on $\r$ is $\n\leq\r+3$, so
$\n-2\leq\r+1$, a contradiction.  Therefore we do in fact know that
$0\neq\prod_{\i=1}^{\n-2}\uu\i\rp\vv\i$, as desired.  This shows $\B_\n(\CEA)$ is a
cylindric basis for $\CEA$ and completes the proof of part~\ref{li}. Parts~\ref{lii}
and~\ref{liii} follow from part~\ref{li} by~\cite[Theorem~10]{MR1011183}.

For part~\ref{liv}, assume to the contrary that $\Ca{\B_\n(\CEA)}\subseteq\Nr_\n\gc\D$
for some $\gc\D\in\CA{\r+4}$.  We get a contradiction by finding a subalgebra $\gc\F$ of
$\Ca{\B_\n(\CEA)}$ which is not in $\Sub\nr_\n\CA{\r+4}$.  From Theorem~\ref{th5} with
$\p=\r+3$ we get
\begin{equation}\label{complexcylalg}
	\Ca{\B_3(\gc\F)}\notin\Sub\mathsf{Nr}_3\CA{\p+1}=\Sub\mathsf{Nr}_3\CA{\r+4}
\end{equation}
For this we choose an arbitrary finite parameter $\N\in\omega$ and make it big enough.
For this fixed $\N$ there is a finite subalgebra $\gc\F$ of $\CEA$ whose atoms are
$\id$, $\pra\i\j$ and $\jn{\cc\i}\N$ for $1\leq\i\leq\r+2$ and $\j<\N$.  This finite
subalgebra $\gc\F$ has a subalgebra isomorphic to $\ep$, whose atoms are $\id$ and
$\jn{\cc\i}0$ for $1\leq\i\leq\r+2$.  By Theorem~\ref{th5} in the next section, a finite
extension of $\ep$ with enough atoms satisfies~\eqref{complexcylalg}. By choosing $\N$
large enough, the extension $\gc\F$ of $\ep$ has enough atoms.
\endproof
\let\germ=\gc

\providecommand{\bysame}{\leavevmode\hbox to3em{\hrulefill}\thinspace}
\providecommand{\MR}{\relax\ifhmode\unskip\space\fi MR }
\providecommand{\MRhref}[2]{%
  \href{http://www.ams.org/mathscinet-getitem?mr=#1}{#2}
}
\providecommand{\href}[2]{#2}

\appendix

\section{\bf A representation result}
The next theorem was invoked in the proof of Theorem~\ref{th2}.
\begin{theorem}\label{th3}
Assume $\gc\A$ is an atomic symmetric integral relation algebra containing a flexible
trio. Then $\gc\A$ has the 1-point extension property and $\gc\A\in\RRA$.
\end{theorem}
\proof
Once it has been shown that $\gc\A$ has the 1-point extension property, it follows from
some additional observations about the behavior of identity elements that the set
$\B_\k(\gc\A)$ of basic $\k$-by-$\k$ matrices of atoms of $\gc\A$ is a relational basis
for $\gc\A$ whenever $\k\geq3$ (see Def.~\ref{def-basis}).  Then $\gc\A\in\RA_\k$ for
all $\k\geq3$ because $\gc\A$ is atomic and has a $\k$-dimensional relational basis,
hence $\gc\A\in\bigcap_{\k\geq3}\RA_\k=\RRA$. In fact, when $\gc\A$ is atomic (as in the
proof of Theorem~\ref{th1}) it is easy to prove directly from the 1-point extension
property that $\gc\A$ has a complete representation on $\omega$.

We show next that there is a function $\f$ such that for any diversity atoms $\x$ and
$\y$, we have $\f(\x,\y)\in\{\a,\b,\c\}$ and
\begin{equation}\label{good}
\x\rp\f(\x,\y)=\y\rp\f(\x,\y)\geq\di
\end{equation}
For an arbitrary diversity atom $\x$, consider the set
\begin{equation}\label{defZ}
\Z_\x:=\{\z:\x\rp\z\geq\di,\,\di\geq\z\in\at{\gc\A}\}.
\end{equation}
If $\x\in\{\a,\b,\c\}$ then $\a,\b,\c\in\Z_\x$ by \eqref{m1}.  If $\x\notin\{\a,\b,\c\}$
then by~\eqref{m}, $\{\a,\b,\c\}\cap\Z_\x$ has at least two elements.  Consequently, if
$\y$ is another, possibly different, diversity atom of $\gc\A$, then, since $\Z_\x$ and
$\Z_\y$ are subsets of the 3-element set $\{\a,\b,\c\}$ and they each contain at least
two elements, they must intersect. We choose a value in the intersection as
$\f(\x,\y)$. There are several ways to do this. We pick this one---
\begin{enumerate}
\item[\bf(A)]	if $\a\in\Z_\x\cap\Z_y$ then $\f(\x,\y)=\a$,
\item[\bf(B)]	if $\a\notin\Z_\x\cap\Z_y$ and $\b\in\Z_\x\cap\Z_y$ then $\f(\x,\y)=\b$,
\item[\bf(C)]	if $\a\notin\Z_\x\cap\Z_y$ and $\b\notin\Z_\x\cap\Z_y$ then $\f(\x,\y)=\c$.
\end{enumerate}
It is obvious from \eqref{defZ} that \eqref{good} holds in the first two cases. We need
to show \eqref{good} also holds in the third case (C), \ie, that under the assumptions
$\a\notin\Z_\x\cap\Z_y$ and $\b\notin\Z_\x\cap\Z_y$ we have $\c\in\Z_\x\cap\Z_\y$.  But
if $\c\notin\Z_\x\cap\Z_\y$ then we would conclude that $\Z_\x\cap\Z_\y$ is empty, since
it is a subset of $\{\a,\b,\c\}$ that excludes each of $\a$, $\b$, and $\c$ by our
assumptions, contrary to the observations made above.

For the 1-point extension property, assume $\k\in\omega$, $\mo\in\B_\k(\gc\A)$,
$\mo_{\l,\m}$ satisfies the identity condition, $\x,\y$ are diversity atoms of $\gc\A$,
and $\mo_{\i,\j}\leq\x\rp\y$ for some fixed $\i,\j<\k$. We will prove that $\mo$ has a
1-point extension $\mo'\in\B_{\k+1}(\gc\A)$ such that $\mo\subseteq\mo'$,
$\mo'_{\i,\k}=\x=\mo'_{\k,\i}$, $\mo'_{\k,\j}=\y=\mo'_{\j,\k}$, and if $\k>\l\neq\i,\j$
then
\begin{equation*}
	\mo'_{\l,\k}=\mo'_{\k,\l}=\f(\x,\y)
\end{equation*}
Note that by definitions and \eqref{good} we have
\begin{equation}\label{key}
	\x\rp\mo'_{\k,\l}=\di=\y\rp\mo'_{\k,\l}.
\end{equation}
Having chosen $\mo'_{\k,\l}$ to be either $\a$ or $\b$ or $\c$, we must check for each
$\l<\k$ whether the first two crucial cycle equations below hold, and finally whether
the third equation holds for those points $\l,\m<\k$ where $\l\neq\m$ and
$\{\l,\m\}\cap\{\i,\j\}=\emptyset$.
\begin{align*}
	\mo_{\i,\l}&\leq\mo'_{\i,\k}\rp\mo'_{\k,\l}
&&\text{\ie, \ $[\mo'_{\i,\k},\mo'_{\k,\l},\mo_{\i,\l}]$ is a cycle}
\\	\mo_{\j,\l}&\leq\mo'_{\j,\k}\rp\mo'_{\k,\l}
&&\text{\ie, \ $[\mo'_{\j,\k},\mo'_{\k,\l},\mo_{\j,\l}]$ is a cycle}
\\	\mo_{\l,\m}&\leq\mo'_{\l,\k}\rp\mo'_{\k,\m}
&&\text{\ie, \ $[\mo'_{\l,\k},\mo'_{\k,\m},\mo_{\l,\m}]$ is a cycle}
\end{align*}
The first two equations hold by \eqref{key} (their right sides are $\di$). For the third
equation, first note that $\mo'_{\l,\k}=\mo'_{\k,\m}$ because the value depends only on
$\x,\y$, not on $\l$ or $\m$.  The right side of the third equation is therefore
$\a\rp\a$ or $\b\rp\b$ or $\c\rp\c$, but $\a\rp\a=\b\rp\b=\c\rp\c=1$, so the third
equation holds.
\endproof

\section{\bf A non-representation result}
\begin{theorem}\label{th5}
Assume 
\begin{enumerate}
\item	$\gc\E\subseteq\gc\A\in\NA$,
\item	$\gc\E$ is finite and symmetric, $\id\in\at{\gc\E}$, and $\gc\E$ has $\p>3$ atoms,
\item	$\gc\E$ has no 1-cycles: $\u\rp\u\bp\u=0$ if $\di\geq\u\in\at{\gc\E}$, 
\item	$\gc\A$ is finite and symmetric, $\id\in\at{\gc\A}$, and some diversity atom of
	$\gc\E$ is the join of at least $\p^{\p-1}$ atoms of $\gc\A$.
\end{enumerate}
Then
\begin{enumerate}
\item	$\gc\A\notin\Sub\mathsf{Ra}\CA{\p+1}$,
\item	$\Ca{\B_3(\gc\A)}\notin\Sub\mathsf{Nr}_3\CA{\p+1}$.
\end{enumerate}
\end{theorem}
\proof[Proof of (i)]
Let the atoms of $\gc\E$ be $\id=\aa0$, $\aa1$, \dots, $\aa{\p-1}$, where $\p\geq3$ and
$\aa1$ is a diversity atom of $\gc\E$ which is the join of at least $\p^{\p-1}$ atoms
of $\gc\A$. Let $\cov\x$ be the atom of $\gc\E$ containing $\x\in\at{\gc\A}$.  We refer
to $\cov\x$ as the ``color'' of $\x$ (or ``cover'', as in the definition of splitting).

Assume, for the sake of obtaining a contradiction, that $\gc\A\subseteq\Ra{\gc\D}$ for
some $\gc\D\in\CA{\p+1}$.  All the elements of $\gc\A$, in particular all the atoms, are
2-dimensional elements of $\gc\D$, \ie,
\begin{equation}\label{nr2}
	\at{\gc\A}\subseteq{Nr}_2\gc\D.
\end{equation}
If $1<\q\leq\p+1$ and $\x\in\D$, we say $\x$ is~\concept{$\q$-color-ordered} if
$\cov\u=\cov\v$ whenever $\u,\v\in\at{\gc\A}$, $0\leq\i<\j<\k<\q$, and
$\x\leq\sub0\i\sub1\j\u\bp\sub0\i\sub1\k\v$.

The element $\x\in\D$ is~\concept{$\q$-covered} if there are atoms
$\uu{\i,\j}\in\at{\gc\A}$ for $0\leq\i<\j<\q$ such that
$\x\leq\prod_{0\leq\i<\j<\q}\sub0\i\sub1\j\uu{\i,\j}$, in which case the atoms
$\uu{\i,\j}$ are said to be a~\concept{$\q$-covering} of $\x$.

The atoms in a $\q$-covering of a non-zero $\x\in\D$ are unique, for if there are
further atoms $\vv{\i,\j}\in\at{\gc\A}$, $0\leq\i<\j<\q$, such that
$\x\leq\prod_{0\leq\i<\j<\q}\sub0\i\sub1\j\vv{\i,\j}$, then, since substitution is a
complete Boolean endomorphism by~\cite[1.5.3]{MR0314620}, we have
\begin{align*}
0\neq\x&\leq
	\prod_{0\leq\i<\j<\q}\sub0\i\sub1\j\uu{\i,\j}
	\bp\prod_{0\leq\i<\j<\q}\sub0\i\sub1\j\vv{\i,\j}
\\&=	\prod_{0\leq\i<\j<\q}\sub0\i\sub1\j(\uu{\i,\j}\bp\vv{\i,\j}),
\end{align*}
but if $\uu{\i,\j}\neq\vv{\i,\j}$ then, since distinct atoms are disjoint, a zero occurs
with a contradiction ensuing. Thus $\uu{\i,\j}=\vv{\i,\j}$ whenever $0\leq\i<\j<\q$.

We will construct by induction for each dimension from $\q=2$ up to $\q=\p+1$ a set
$\S_\q\subseteq{Nr}_\q\gc\D$ such that
\begin{enumerate}
\item\label{p1}	$\S_\q$ has at least $\p^{\p+1-\q}$ elements.
\item\label{p2}	Every $\x\in\S_\q$ is $\q$-covered, $\q$-color-ordered, and non-zero,
		and $\x\leq\sub1\j\aa1$ for $0<\j<\q$.
\item\label{p3}	$\cyl{\q-1}\x=\cyl{\q-1}\y$ if $\x,\y\in\S_\q$.
\item\label{p4}	$\cov\u=\cov\v$ if $\u,\v\in\at{\gc\A}$, $\x,\y\in\S_\q$,
		$\x\leq\sub0{\q-2}\sub1{\q-1}\u$, and $\y\leq\sub0{\q-2}\sub1{\q-1}\v$.
\item\label{p5}	$\u\neq\v$ if $\u,\v\in\at{\gc\A}$, $\x,\y\in\S_\q$,
		$\x\leq\sub1{\q-1}\u$, $\y\leq\sub1{\q-1}\v$, and $\x\neq\y$.
\item\label{p6}	$\u\neq\v$ if $0<\j<\k<\q$, $\u,\v\in\at{\gc\A}$, $\x\in\S_\q$, and
		$\x\leq\sub1\j\u\bp\sub1\k\v$.
\end{enumerate}
Let $\S_2=\{\x:\aa1\geq\x\in\at{\gc\A}\}$.  

Note that $\S_2\subseteq{Nr}_2\gc\D$ by~\eqref{nr2}.  Obviously $\S_2$ has
property~\ref{p1} since there are at least $\p^{\p-1}$ atoms below $\aa1$.  Let
$\x\in\S_2$. Then $\x$ is $2$-covered by itself (take $\uu{0,1}=\x$), $\x$ is
$2$-color-ordered because the hypotheses in the definition of color-ordered are never
met ($\q=2$ is too small), and $\x$ is not zero because it is an atom of $\gc\A$. For
the last part of property~\ref{p2}, note that if $0\leq\i<\j<\q=2$ then $\j=1$, and
$\x\leq\aa1$ by the definition of $\S_2$, so $\x\leq\aa1=\sub11\aa1=\sub1\j\aa1$.
Therefore $\S_2$ has property~\ref{p2}.  Since $\gc\A$ is integral and $\x\in\S_2$ is
non-zero, we have $\x\rp1=1$, so
\begin{align*}
 1=\x\rp1
	&=\cyl2(\sub12\x\bp\sub021)	&&\text{definition of $\rp$ in $\Ra{\gc\D}$}
\\	&=\cyl2\sub12\x			&&\text{\cite[1.5.3]{MR0314620}}
\\	&=\cyl1\sub21\x			&&\text{\cite[1.5.9(i)]{MR0314620}}
\\	&=\cyl1\x			&&\text{\cite[1.5.8(i)]{MR0314620}, $\cyl2\x=\x$}
\end{align*}
It follows that property~\ref{p3} holds for $\S_2$.  For property~\ref{p4}, note that
since $\q=2$, $\sub0{\q-2}\sub1{\q-1}$ is the identity mapping, hence the hypotheses are
$\u,\v\in\at{\gc\A}$, $\x,\y\in\S_2$, $\x\leq\u$, and $\y\leq\v$, which imply $\x=\u$
and $\y=\v$ since $\u,\v,\x,\y$ are atoms.  We wish to show $\cov\u=\cov\v$, \ie,
$\cov\x=\cov\y$, but this is true by the definition of $\S_2$.  Since $\q=2$, the
substitution $\sub1{\q-1}$ is the identity mapping, hence the hypotheses of
property~\ref{p5} are $\u,\v\in\at{\gc\A}$, $\x,\y\in\S_\q$, $\x\leq\u$, $\y\leq\v$, and
$\x\neq\y$. But these hypotheses imply $\u=\x\neq\y=\v$, so the conclusion holds
trivially. Thus $\S_2$ has property~\ref{p5}.  Finally, $\S_2$ has property~\ref{p6}
because the hypotheses cannot hold when $\q=2$.

Suppose we have a set $\S_\q$ such that $\q\geq2$ and \ref{p1}--\ref{p6}.  Choose an
arbitrary but fixed $\w\in\S_\q$, and let $\sq:=\S_\q\minus\{\w\}$.  We will obtain a
function $\h$ that sends every $\x\in\sq$ to a $(\q+1)$-dimensional element
$\h(\x)\in{Nr}_{\q+1}\gc\D$, and will choose $\S_{\q+1}$ to be a subset of the range of
$\h$.

For every $\x\in\sq$, we have
\begin{align*}
0	&\neq\w			&&\text{\ref{p2}}
\\	&=\w\bp\cyl{\q-1}\w	&&\text{\cite[1.1.1(C$_2$)]{MR0314620}}
\\	&=\w\bp\cyl{\q-1}\x	&&\text{\ref{p3}}
\\	&=\w\bp\cyl{\q-1}\sub\q{\q-1}\x
			&&\text{\cite[1.5.8(i)]{MR0314620}, $\cyl\q\x=\x$}
\\	&=\w\bp\cyl\q\sub{\q-1}\q\x
			&&\text{\cite[1.5.9(i)]{MR0314620}}
\\	&=\cyl\q(\w\bp\sub{\q-1}\q\x)	
			&&\text{\cite[1.1.1(C$_3$)]{MR0314620}, $\cyl\q\w=\w$}
\\	&=\cyl\q\(\w\bp\sub{\q-1}\q\x\bp\sub0{\q-1}\sub1\q\(1\)\)
			&&\text{\cite[1.5.3]{MR0314620}}
\\	&=\cyl\q\(\w\bp\sub{\q-1}\q\x\bp\sub0{\q-1}\sub1\q\(\sum_{\y\in\at{\gc\A}}\y\)\)
			&&\text{$\gc\A$ is finite}
\\	&=\sum_{\y\in\at{\gc\A}}\cyl\q\(\w\bp\sub{\q-1}\q\x\bp\sub0{\q-1}\sub1\q(\y)\)
			&&\text{\cite[1.5.3, 1.2.6]{MR0314620}}
\end{align*}
The distributive law holds in all Boolean algebras whenever all the joins and meets
involved are finite, so
\begin{align*}
0\neq\w
	&=\prod_{\x\in\sq}
	\(\sum_{\y\in\at{\gc\A}}\cyl\q\(\w\bp
	\sub{\q-1}\q\x\bp\sub0{\q-1}\sub1\q(\y)\)\)
\\	&=\sum_{\f:\sq\to\at{\gc\A}}\(\prod_{\x\in\sq}
	\cyl\q\(\w\bp\sub{\q-1}\q\x\bp\sub0{\q-1}\sub1\q(\f(\x))\)\).
\end{align*}
Consequently there must be some function $\f:\sq\to\at{\gc\A}$ such that
\begin{equation}\label{key2}
0\neq\prod_{\x\in\sq}\cyl\q\(\w\bp\sub{\q-1}\q\x\bp\sub0{\q-1}\sub1\q(\f(\x))\).
\end{equation}
Let $\f$ be such a function. From our chosen $\f$ we define additional functions
$\g,\h:\sq\to\D$ and an element $\z\in\D$ as follows.
\begin{align}
	\g(\x)	&=\w\bp\sub{\q-1}\q\x\bp\sub0{\q-1}\sub1\q(\f(\x)) 
		&&\text{for all $\x\in\sq$}\label{8}
\\	\z	&=\prod_{\x\in\sq}\cyl\q\(\g(\x)\)\label{9}
\\	\h(\x)	&=\g(\x)\bp\z
		&&\text{for all $\x\in\sq$}\label{10}
\end{align}
Let $\R=\{\h(\x):\x\in\sq\}$. We will show that $\R$ itself has properties~\ref{p2},
\ref{p3}, ~\ref{p5}, and~\ref{p6}. Consequently every subset of $\R$ also has these
properties.  We will partition $\R$ into disjoint subsets that have property~\ref{p4}
and prove that at least one of them must be large enough to also have property~\ref{p1}.
We take $\S_{\q+1}$ to be any such subset of $\R$.

To see that $\R$ has property~\ref{p3}, we observe that $\cyl\q\h(\x)=\cyl\q\h(\y)$ for
all $\x,\y\in\sq$, because
\begin{align*}
	\cyl\q\h(\x)
  &=	\cyl\q(\g(\x)\bp\z)	&&\text{\eqref{10}}
\\&=	\cyl\q(\g(\x))\bp\z	&&\text{\cite[1.1.1(C$_3$)]{MR0314620}, $\cyl\q\z=\z$}
\\&=	\z			&&\text{\eqref{9}}
\end{align*}
It follows that $\h(\x)\neq0$ for every $\x\in\S_\q$, since $\z\neq0$
by~\eqref{key2}. This is part of property~\ref{p2}. For the last part of
property~\ref{p2}, we want to show $\h(\x)\leq\sub1\j(\aa1)$ whenever $0<\j<\q+1$ and
$\x\in\S_2$.  We have $\h(\x)\leq\g(\x)\leq\w\bp\sub{\q-1}\q\x$ by
definitions~\eqref{10} and~\eqref{8}, so there are two cases. First, assume $0<\j<\q$.
In this case we note that from $\w\in\S_\q$ and~\ref{p2} for $\S_\q$ we get
$\w\leq\sub1\j\aa1$, so $\h(\x)\leq\sub1\j\aa1$. Suppose $\j=\q$.  In this case we have
$\x\leq\sub1\j\aa1$ for $0<\j<\q$ by \ref{p2} for $\S_\q$ since $\x\in\S_\q$.  In
particular, $\x\leq\sub1{\q-1}\aa1$, so
$\h(\x)\leq\sub{\q-1}\q\sub1{\q-1}\aa1=\sub1\q\aa1$.  We get the rest of
property~\ref{p2} by showing $\h(\x)$ is $(\q+1)$-color-ordered and $(\q+1)$-covered for
every $\x\in\sq$.  From $\x\in\sq$ and property~\ref{p2} for $\S_\q$ we know $\x$ is
$\q$-covered, so there are atoms $\x_{\i,\j}\in\at{\gc\A}$ such that
\begin{equation}\label{cov-x}
	\x\leq\prod_{0\leq\i<\j<\q}\sub0\i\sub1\j(\x_{\i,\j}).
\end{equation}
Of course, we also know $\w\in\S_\q$, so there is a $\q$-covering
$\w_{\i,\j}\in\at{\gc\A}$, $0\leq\i<\j<\q$, of $\w$ as well, where
\begin{equation}\label{cov-w}
	\w\leq\prod_{0\leq\i<\j<\q}\sub0\i\sub1\j(\w_{\i,\j}).
\end{equation}
Let
\begin{align}\label{u-def}
\hc{\i,\j}=\begin{cases}
	\w_{\i,\j} 	&\text{if $0\leq\i<\j<\q$}
\\	\x_{\i,\q-1}	&\text{if $0\leq\i<\q-1$ and $\j=\q$}
\\	\f(\x)		&\text{if $\i=\q-1$ and $\j=\q$}
\end{cases}
\end{align}
We shall see that $\hc{\i,\j}$ is a $(\q+1)$-covering of $\h(\x)$.  First, note that
\begin{align}\label{yes}
\sub{\q-1}\q\x&\leq\prod_{0\leq\i<\q-1}\sub0\i\sub1\q(\x_{\i,\q-1})
\end{align}
because if $0\leq\i<\q-1$ then $\x\leq\sub0\i\sub1{\q-1}(\x_{\i,\q-1})$ by
\eqref{cov-x}, so
\begin{align*}
	\sub{\q-1}\q\x
&\leq	\sub{\q-1}\q\sub0\i\sub1{\q-1}(\x_{\i,\q-1})
	&&\text{\cite[1.5.3]{MR0314620}}
\\&=	\sub0\i\sub{\q-1}\q\sub1{\q-1}(\x_{\i,\q-1})
	&&\text{\cite[1.5.10(iii)]{MR0314620}}\notag
\\&=	\sub0\i\sub1\q(\x_{\i,\q-1})
	&&\text{$\cyl{\q-1}\x_{\i,\q-1}=\x_{\i,\q-1}$,
	\cite[1.5.11(i)]{MR0314620}}\notag
\end{align*}
Then we have
\begin{align}\label{cov-h}
	\h(\x)
  &\leq	\g(\x)=\w\bp\sub{\q-1}\q\x\bp\sub0{\q-1}\sub1\q(\f(\x))
	&&\text{\eqref{10}, \eqref{8}}
\\&\leq	\prod_{0\leq\i<\j<\q}\sub0\i\sub1\j(\w_{\i,\j})
	\bp\prod_{0\leq\i<\q-1}\sub0\i\sub1\q(\x_{\i,\q-1})
	\bp\sub0{\q-1}\sub1\q(\f(\x)) 
	&&\text{\eqref{cov-w}, \eqref{yes}}\notag
\\&=	\prod_{0\leq\i<\j<\q}\sub0\i\sub1\j(\hc{\i,\j})
	\bp\prod_{0\leq\i<\q-1}\sub0\i\sub1\q(\hc{\i,\q})
	\bp\sub0{\q-1}\sub1\q(\hc{\q-1,\q})
	&&\text{\eqref{u-def}}\notag
\\&=	\prod_{0\leq\i<\j<\q+1}\sub0\i\sub1\j(\hc{\i,\j})
\notag
\end{align}
so $\h(\x)$ is $(\q+1)$-covered.  

To show $\h(\x)$ is $(\q+1)$-color-ordered, we assume $0\leq\i<\j<\k<\q+1$ and must show
$\cov{\hc{\i,\j}}=\cov{\hc{\i,\k}}$.  If $\i<\j<\k<\q$ then the first case
in~\eqref{u-def} applies to both $\hc{\i,\j}$ and $\hc{\i,\k}$, hence
$\hc{\i,\j}=\w_{\i,\j}$ and $\w_{\i,\k}=\hc{\i,\k}$, but
$\cov{\w_{\i,\j}}=\cov{\w_{\i,\k}}$ because $\w$ is color-ordered, so we have
$\cov{\hc{\i,\j}}=\cov{\hc{\i,\k}}$. We may therefore assume $\k=\q$.

We need to observe before going on that if $\q>2$, then
\begin{align*}
\w\leq\cyl{\q-1}\w
  &= \cyl{\q-1}\x	&&\text{property~\ref{p3} of $\S_\q$}
\\&=	\cyl{\q-1}\(
	\prod_{0\leq\i<\j<\q-1}\sub0\i\sub1\j(\x_{\i,\j})\)
	&&\text{\cite[1.2.6]{MR0314620}, \eqref{cov-x}}
\\&=	\prod_{0\leq\i<\j<\q-1}\sub0\i\sub1\j(\x_{\i,\j})
	&&\text{$\cyl{\q-1}\x_{\i,\j}=\x_{\i,\j}$, $\q-1\geq2$}
\end{align*}
By the uniqueness of coverings this tells us that 
\begin{equation}\label{reason}
\w_{\i,\j}=\x_{\i,\j}\text{ if }0\leq\i<\j<\q-1.
\end{equation}
If, in addition to $\k=\q$, we have $\i<\j<\q-1$, then $\q>2$ and the first and second
cases of~\eqref{u-def} apply, so we have $\hc{\i,\j}=\w_{\i,\j}$ and
$\hc{\i,\k}=\hc{\i,\q}=\x_{\i,\q-1}$.  But $\w_{\i,\j}=\x_{\i,\j}$ by~\eqref{reason}.
Also, $\x$ is color-ordered, so $\cov{\x_{\i,\j}}=\cov{\x_{\i,\q-1}}$, which is
equivalent to $\cov{\hc{\i,\j}}=\cov{\hc{\i,\k}}$ by the previous equations.

The final case is that $\i<\j=\q-1$ and $\k=\q$. The possibilities for $\i$ divide into
two sub-cases, $\i$ is smaller than $\q-2$, and $\i$ is equal to $\q-2$. If
$0\leq\i<\q-2$ then $\i<\q-2<\q-1$, so $\cov{\hc{\i,\j}} =\cov{\w_{\i,\q-1}}
=\cov{\w_{\i,\q-2}}$ since $\w$ is color-ordered by property~\ref{p2} of $\S_\q$, and
$\cov{\x_{\i,\q-2}}=\cov{\x_{\i,\q-1}}$ since $\x$ is color-ordered, but
$\w_{\i,\q-2}=\x_{\i,\q-2}$ by~\eqref{reason}, so
\begin{align*}
	\cov{\hc{\i,\j}}
  &=	\cov{\w_{\i,\q-1}}	&&\text{$\j=\q-1$, \eqref{u-def}}
\\&=	\cov{\w_{\i,\q-2}}	&&\text{$\w$ is color-ordered}
\\&=	\cov{\x_{\i,\q-2}}	&&\text{\eqref{reason}}
\\&=	\cov{\x_{\i,\q-1}}	&&\text{$\x$ is color-ordered}
\\&=	\cov{\hc{\i,\k}}	&&\text{$\q=\k$, \eqref{u-def}}
\end{align*}
We are reduced to assuming $\i=\q-2$, hence
\begin{equation*}
\cov{\hc{\i,\j}}=\cov{\w_{\q-2,\q-1}}=\cov{\x_{\q-2,\q-1}}=\cov{\hc{\i,\k}}
\end{equation*}
by~\eqref{reason} and the third case in \eqref{u-def}.  

We have shown that every $\h(\x)$ constructed from some $\x\in\sq$ is non-zero,
$(\q+1)$-covered, and $(\q+1)$-color-ordered. Thus $\R$ and all its subsets has
property~\ref{p2}.

To prove property~\ref{p5} for $\R$ (and its subsets), we assume $\x,\y\in\sq$,
$\h(\x)\neq\h(\y)$, $\u,\v\in\at{\gc\A}$, $\h(\x)\leq\sub1\q\u$, $\h(\y)\leq\sub1\q\v$. We
must show $\u\neq\v$.  If we have a $\q$-covering of $\x$ as in~\eqref{cov-x}, then
by~\eqref{cov-h} we get $\u=\x_{0,\q-1}$ from $\h(\x)\leq\sub1\q\u$, and, similarly,
$\v=\y_{0,\q-1}$ from $\h(\y)\leq\sub1\q\v$ for some $\q$-covering $\y_{\i,\j}$ of $\y$.
Hence $\x\leq\sub1{\q-1}(\x_{0,\q-1})$ and $\y\leq\sub1{\q-1}(\y_{0,\q-1})$, so, by
property~\ref{p5} for $\sq$, we know $\x_{0,\q-1}\neq\y_{0,\q-1}$, \ie, $\u\neq\v$, as
desired.

To prove property~\ref{p6} for $\R$ (and its subsets), we assume $0<\j<\k<\q+1$,
$\u,\v\in\at{\gc\A}$, $\x\in\sq$, and $\h(\x)\leq\sub1\j\u\bp\sub1\k\v$.  If $\k<\q$,
then $\u=\hc{0,\j}=\w_{0,\j}$ and $\v=\hc{0,\k}=\w_{0,\k}$ by~\eqref{cov-h}
and~\eqref{u-def}, but $\w\in\S_\q$, so by property~\ref{p6} for $\S_\q$, we have
$\w_{0,\j}\neq\w_{0,\k}$, hence $\u\neq\v$.  Suppose that $\k=\q$.  In this case,
by~\eqref{cov-h} and~\eqref{u-def}, we again have $\u=\hc{0,\j}=\w_{0,\j}$ but this time
$\v=\hc{0,\q}=\x_{0,\q-1}$.  Hence $\w\leq\sub1\j\u$ and $\x\leq\sub1{\q-1}\v$
by~\eqref{cov-x} and~\eqref{cov-w}.  If $\j=\q-1$ we note that $\w\neq\x$ since
$\x\in\sq$, hence $\u\neq\v$ by property~\ref{p5} for $\S_\q$, which gives us
$\hc{0,\j}\neq\hc{0,\q}$, \ie, $\u\neq\v$.  If $\j<\q-1$ then
$\v=\hc{0,\q}=\x_{0,\q-1}\neq\x_{0,\j}$ for $0<\j<\q-1$ by property~\ref{p6} for
$\S_\q$, applied this time to $\x$.  But $\x_{0,\j}=\w_{0,\j}=\hc{0,\j}$
by~\eqref{reason} and~\eqref{u-def}, so again we have $\hc{0,\q}\neq\hc{0,\j}$.

We have proved $\R$ has properties~\ref{p2}, \ref{p3}, ~\ref{p5}, and~\ref{p6}, and wish
to show that $\h$ is one-to-one on $\sq$. Assume $\x,\y\in\sq$ and $\x\neq\y$. We want
to show $\h(x)\neq\h(\y)$.  By property~\ref{p2} for $\sq$, $\x$ and $\y$ have
$\q$-coverings that include atoms $\x_{0,\q-1},\y_{0,\q-1}\in\at{\gc\A}$ satisfying
$\x\leq\sub1{\q-1}(\x_{0,\q-1})$ and $\y\leq\sub1{\q-1}(\y_{0,\q-1})$.  By
\eqref{u-def} and~\eqref{cov-h} these last two equations imply
$\h(\x)\leq\sub1\q(\x_{0,\q-1})$ and $\h(\y)\leq\sub1\q(\y_{0,\q-1})$.  From $\x\neq\y$
we conclude by property~\ref{p5} for $\sq$ that $\x_{0,\q-1}\neq\y_{0,\q-1}$, which
implies, by property~\ref{p5} for $\R$, that $\h(x)\neq\h(\y)$, as desired.

Now we want to choose a subset $\S_{\q+1}$ of $\R$ with property~\ref{p4} that contains
at least $\p^{\p+1-(\q+1)}$ elements.  We partition $\R$ and let $\S_{\q+1}$ be the
largest piece.  Recall from \eqref{cov-h} that
$\h(\x)\leq\sub0{\q-2}\sub1{\q-1}(\f(\x))$ for every $\x\in\sq$, and $\f(\x)$ has color
$\cov{\f(\x)}\in\at{\gc\E}$.  For every color $\aa\i$ we get a piece of $\R$, namely
\begin{equation*}
	\R_\i:=\{\h(\x):\x\in\sq,\,\cov{\f(\x)}=\aa\i\}.
\end{equation*}
Note that $\R$ is the disjoint union of the pieces, the number of pieces is $\p$, and
$\R$ has at least $\p^{\p+1-\q}$ elements because $\h$ is one-to-one and $\sq$ has more
than $\p^{\p+1-\q}$ elements.  Consequently some piece has at least
$\p^{\p+1-\q}/\p=p^{\p-\q}$ elements in it, and we let $\S_{\q+1}$ be any such
piece. Thus $\S_{\q+1}$ has property~\ref{p1}. Every piece has property~\ref{p4}, so in
particular $\S_{\q+1}$ has this property. Finally, as a subset of $\R$, $\S_{\q+1}$ has
all the other properties.  This completes the construction of the sets $\S_\q$.

Consider what happens when $\q=\p+1$.  We may choose some $\x\in\S_{\p+1}$ because
$\S_{\p+1}$ has at least one element, by property~\ref{p1}.  Then $\x$ is
$(\p+1)$-covered, $(\p+1)$-color-ordered, and non-zero by property~\ref{p2}. Let $\x$
have $(\p+1)$-covering $\x_{\i,\j}\in\at{\gc\A}$ for $0\leq\i<\j<\p+1$.

Consider the set $\{\cov{\x_{\i,\p}}:0\leq\i<\p\}\subseteq\at{\gc\E}$.  Note that
$\cov{\x_{0,\p}}=\aa1\neq\id$ since $\x_{0,\p}\leq\aa1$ by property~\ref{p2}.  We can
also show $\cov{\x_{\i,\p}}\neq\id$ for $0<\i<\p$ because we have, by the covering of
$\x$, $\x\leq\sub1\i(\x_{0,\i})\bp\sub0\i\sub1\p(\x_{\i,\p})\bp\sub1\p(\x_{0,\p})$ so it
follows by \cite[Lemma~10]{MR1931348} that $[\x_{0,\i},\x_{\i,\p},\x_{0,\p}]$ is a
cycle, \ie, $\x_{0,\i}\rp\x_{\i,\p}\geq\x_{0,\p}$. If $\cov{\x_{\i,\p}}=\id$ then
$\x_{\i,\p}=\id$ and we would get $\x_{0,\i}=\x_{0,\p}$, contradicting
property~\ref{p6}, which says $\x_{0,\i}\neq\x_{0,\p}$ for $0<\i<\p$.  Thus we know
$\cov{\x_{\i,\p}}$ is a diversity atom of $\gc\E$ for $0\leq\i<\p$.

The number of diversity atoms in $\gc\E$ is $\p-1$, but the size of the index set
$\{\i:0\leq\i<\p\}$ is $\p$. Therefore some atom is repeated, \ie, there are
$0\leq\i<\j<\p$ such that $\cov{\x_{\i,\p}}=\cov{\x_{\j,\p}}$.  By the
$(\p+1)$-color-ordering of $\x$, $\cov{\x_{\i,\j}}=\cov{\x_{\i,\p}}$.  Let
$\u=\cov{\x_{\i,\j}}=\cov{\x_{\i,\p}}=\cov{\x_{\j,\p}}$. We proved above that
$\u\neq\id$. By the covering of $\x$ and property~\ref{p2} we have
$0\neq\x\leq\sub0\i\sub1\j(\x_{\i,\j})\bp\sub0\j\sub1\p(\x_{\j,\p})
\bp\sub0\i\sub1\p(\x_{\i,\p})$, hence by the definition of $\u$
and~\cite[Lemma~10]{MR1931348} we have $0\neq\u\rp\u\bp\u$. Since $\u\neq\id$, this
contradicts the assumption that $\gc\E$ has no such diversity atom as the $\u$ we have
found.
\endproof
\proof[Proof of (ii)]
Assume to the contrary that $\Ca{\B_3(\gc\A)}\in\Sub\mathsf{Nr}_3\CA{\p+1}$.  Then
\begin{align*}
\gc\A\cong\Ra{\Ca{\B_3(\gc\A)}}	&\in\gc{Ra}^*\Sub\mathsf{Nr}_3\CA{\p+1}
\\				&=\Sub\gc{Ra}^*\mathsf{Nr}_3\CA{\p+1}
				&&\text{\cite[5.3.13]{MR781930}}
\\				&=\Sub\gc{Ra}^*\CA{\p+1}
				&&\text{\cite{MR781930}, defs}
\end{align*}
contradicting part~(i) that says $\gc\A\notin\Sub\mathsf{Ra}\CA{\p+1}$.
\endproof


\begin{thebibliography}{10}

\bibitem{MR1052567}
H.~Andr{\'e}ka, R.~D. Maddux, and I.~N{\'e}meti, \emph{Splitting in relation
  algebras}, Proc. Amer. Math. Soc. \textbf{111} (1991), no.~4, 1085--1093.
  \MR{1052567 (91g:03126)}

\bibitem{MR2387933}
Hajnal Andr{\'e}ka, Istv{\'a}n N{\'e}meti, and Tarek Sayed~Ahmed,
  \emph{Omitting types for finite variable fragments and complete
  representations of algebras}, J. Symbolic Logic \textbf{73} (2008), no.~1,
  65--89. \MR{2387933 (2008m:03130)}

\bibitem{MR0357172}
Claude Berge, \emph{Graphs and hypergraphs}, North-Holland Publishing Co.,
  Amsterdam, 1973, Translated from the French by Edward Minieka, North-Holland
  Mathematical Library, Vol. 6. \MR{0357172 (50 \#9640)}

\bibitem{MR734546}
Stephen~D. Comer, \emph{Color schemes forbidding monochrome triangles},
  Proceedings of the fourteenth {S}outheastern conference on combinatorics,
  graph theory and computing ({B}oca {R}aton, {F}la., 1983), vol.~39, 1983,
  pp.~231--236. \MR{734546 (85f:05052)}

\bibitem{Comer1984}
\bysame, \emph{Combinatorial aspects of relations}, Algebra Universalis
  \textbf{18} (1984), no.~1, 77--94.

\bibitem{MR1608984}
M.~F. Frias and R.~D. Maddux, \emph{Non-embeddable simple relation algebras},
  Algebra Universalis \textbf{38} (1997), no.~2, 115--135. \MR{1608984
  (99a:03063)}

\bibitem{MR0314620}
Leon Henkin, J.~Donald Monk, and Alfred Tarski, \emph{Cylindric algebras.
  {P}art {I}. {W}ith an introductory chapter: {G}eneral theory of algebras},
  North-Holland Publishing Co., Amsterdam, 1971, Studies in Logic and the
  Foundations of Mathematics, Vol. 64. \MR{0314620 (47 \#3171)}

\bibitem{MR781930}
\bysame, \emph{Cylindric algebras. {P}art {II}}, Studies in Logic and the
  Foundations of Mathematics, vol. 115, North-Holland Publishing Co.,
  Amsterdam, 1985. \MR{781930 (86m:03095b)}

\bibitem{MR1330986}
Robin Hirsch, \emph{Completely representable relation algebras}, Bull. IGPL
  \textbf{3} (1995), no.~1, 77--91. \MR{1330986 (96d:03082)}

\bibitem{MR1472125}
Robin Hirsch and Ian Hodkinson, \emph{Complete representations in algebraic
  logic}, J. Symbolic Logic \textbf{62} (1997), no.~3, 816--847. \MR{1472125
  (98m:03123)}

\bibitem{MR1935083}
\bysame, \emph{Relation {A}lgebras by {G}ames}, Studies in Logic and the
  Foundations of Mathematics, vol. 147, North-Holland Publishing Co.,
  Amsterdam, 2002, With a foreword by Wilfrid Hodges. \MR{1935083
  (2003m:03001)}

\bibitem{MR1887031}
\bysame, \emph{Strongly representable atom structures of relation algebras},
  Proc. Amer. Math. Soc. \textbf{130} (2002), no.~6, 1819--1831. \MR{1887031
  (2002k:03113)}

\bibitem{MR2548463}
\bysame, \emph{Strongly representable atom structures of cylindric algebras},
  J. Symbolic Logic \textbf{74} (2009), no.~3, 811--828. \MR{2548463
  (2011c:03147)}

\bibitem{MR1931348}
Robin Hirsch, Ian Hodkinson, and Roger~D. Maddux, \emph{Provability with
  finitely many variables}, Bull. Symbolic Logic \textbf{8} (2002), no.~3,
  348--379. \MR{1931348 (2003m:03013)}

\bibitem{MR1490103}
Ian Hodkinson, \emph{Atom structures of cylindric algebras and relation
  algebras}, Ann. Pure Appl. Logic \textbf{89} (1997), no.~2-3, 117--148.
  \MR{1490103 (99c:03103)}

\bibitem{MR2156722}
Ian Hodkinson and Yde Venema, \emph{Canonical varieties with no canonical
  axiomatisation}, Trans. Amer. Math. Soc. \textbf{357} (2005), no.~11,
  4579--4605. \MR{2156722 (2006e:03106)}

\bibitem{MR2605423}
Mohamed Khaled and Tarek Sayed~Ahmed, \emph{Classes of algebras that are not
  closed under completions}, Bull. Sect. Logic Univ. \L\'od\'z \textbf{38}
  (2009), no.~1-2, 29--43. \MR{2605423 (2011a:03071)}

\bibitem{MR2507303}
\bysame, \emph{On complete representations of algebras of logic}, Log. J. IGPL
  \textbf{17} (2009), no.~3, 267--272. \MR{2507303 (2011a:03072)}

\bibitem{MR0037278}
Roger~C. Lyndon, \emph{The representation of relational algebras}, Ann. of
  Math. (2) \textbf{51} (1950), 707--729. \MR{0037278 (12,237a)}

\bibitem{MR662049}
Roger Maddux, \emph{Some varieties containing relation algebras}, Trans. Amer.
  Math. Soc. \textbf{272} (1982), no.~2, 501--526. \MR{662049 (84a:03079)}

\bibitem{Maddux1978}
Roger~D. Maddux, \emph{Topics in relation algebras}, Ph.D. thesis, University
  of California, Berkeley, 1978, pp.~iii+241.

\bibitem{Maddux1985}
\bysame, \emph{Finite integral relation algebras}, Universal Algebra and
  Lattice Theory (Charleston, S.C., 1984), Springer, Berlin, 1985,
  pp.~175--197.

\bibitem{MR1011183}
\bysame, \emph{Nonfinite axiomatizability results for cylindric and relation
  algebras}, J. Symbolic Logic \textbf{54} (1989), no.~3, 951--974. \MR{1011183
  (90f:03099)}

\bibitem{MR2269199}
\bysame, \emph{Relation {A}lgebras}, Studies in Logic and the Foundations of
  Mathematics, vol. 150, Elsevier B. V., Amsterdam, 2006. \MR{2269199
  (2007j:03096)}

\bibitem{MR0256861}
J.~Donald Monk, \emph{Nonfinitizability of classes of representable cylindric
  algebras}, J. Symbolic Logic \textbf{34} (1969), 331--343. \MR{0256861 (41
  \#1517)}

\bibitem{MR0277369}
\bysame, \emph{Completions of {B}oolean algebras with operators}, Math. Nachr.
  \textbf{46} (1970), 47--55. \MR{0277369 (43 \#3102)}

\bibitem{MR2357188}
Tarek Sayed~Ahmed, \emph{Neat embedding is not sufficient for complete
  representability}, Bull. Sect. Logic Univ. \L\'od\'z \textbf{36} (2007),
  no.~1-2, 29--35. \MR{2357188 (2008g:03103)}

\bibitem{MR2443834}
\bysame, \emph{A note on atom structures of relation and cylindric algebras},
  Int. J. Algebra \textbf{2} (2008), no.~9-12, 595--601. \MR{2443834
  (2009g:03101)}

\bibitem{MR2441102}
\bysame, \emph{A note on atom structures of relation and cylindric algebras},
  Bull. Sect. Logic Univ. \L\'od\'z \textbf{37} (2008), no.~1, 29--35.
  \MR{2441102 (2009g:03100)}

\bibitem{MR2417802}
\bysame, \emph{Weakly representable atom structures that are not strongly
  representable, with an application to first order logic}, MLQ Math. Log. Q.
  \textbf{54} (2008), no.~3, 294--306. \MR{2417802 (2009d:03161)}

\bibitem{MR2519240}
\bysame, \emph{A simple construction of representable relation algebras with
  non-representable completions}, MLQ Math. Log. Q. \textbf{55} (2009), no.~3,
  237--244. \MR{2519240 (2010i:03072)}

\bibitem{MR2453362}
Tarek Sayed~Ahmed and Basim Samir, \emph{The class {$S\germ{Nr}_3{\bf CA}_k$}
  is not closed under completions}, Log. J. IGPL \textbf{16} (2008), no.~5,
  427--429. \MR{2453362}

\bibitem{MR0066303}
Alfred Tarski, \emph{Contributions to the theory of models. {III}}, Nederl.
  Akad. Wetensch. Proc. Ser. A. \textbf{58} (1955), 56--64 = Indagationes Math.
  17, 56--64 (1955). \MR{0066303 (16,554h)}

\end{thebibliography}
\end{document}